\documentclass[11pt]{article}

\usepackage{cite}
\usepackage{subfigure}
\usepackage{graphicx, color, graphpap}
\usepackage{amsmath,amssymb,amsthm}
\usepackage{ifthen}
\usepackage{mathrsfs}
\usepackage{mathtools}
\usepackage{enumitem}
\usepackage{lineno}
\usepackage{float}
\usepackage{fancyhdr}
\usepackage[us,12hr]{datetime} 
\usepackage{exscale}
\usepackage{tabularx}
\usepackage{latexsym}
\usepackage{fullpage}       
\usepackage{float}
\usepackage{verbatim}
\usepackage{multirow}
\usepackage{overpic}
\usepackage{comment} 
\usepackage[hidelinks]{hyperref}
\usepackage{blkarray}

\usepackage{tcolorbox}

\usepackage{pgfplots}
\usepackage{marvosym}
\usepackage[noabbrev]{cleveref}
\usepackage{tikz}
\usepackage{pgfplots}
\usepackage{algorithm}
\usepackage{algpseudocode}

\pgfplotsset{compat=1.11}
\usetikzlibrary{arrows, arrows.meta}

\newcommand{\rank}{\text{rank }}
\newcommand{\supp}{\text{supp}}

\newcommand{\tr}{\text{tr}}
\newcommand{\Snp}{\mathbb{S}_{+}^{n}}

\newcommand{\Srp}{\mathbb{S}_{+}^{r}}
\newcommand{\Sn}{\mathbb{S}^{n}}

\newcommand{\msd}{\text{MSD}}
\newcommand{\sd}{\text{SD}}

\newcommand{\Rn}{\R^{n}}
\newcommand{\Rnp}{\R_{+}^{n}}
\newcommand{\Rm}{\R^{m}}

\newcommand{\R}{\mathbb{R}}

\newcommand{\face}{\text{face}}
\newcommand{\spa}{\text{span}}

\newcommand{\conv}{\text{conv}}

\DeclareMathOperator{\ri}{\text{ri}}

\newtheorem{thm}{Theorem}[section]

\newtheorem{prop}[thm]{Proposition}
\newtheorem{lem}[thm]{Lemma}
\newtheorem{cor}[thm]{Corollary}

\newtheorem{assumption}[thm]{Assumption}
\newtheorem{ex}[thm]{Example}

\makeatletter
\let\c@algorithm\c@thm

\makeatother

\crefname{thm}{Theorem}{Theorems}
\Crefname{thm}{Theorem}{Theorems}
\crefname{problem}{Problem}{Theorems}
\Crefname{problem}{Problem}{Theorems}
\Crefname{assump}{Assumption}{Theorems}
\crefname{assump}{Assumption}{Theorems}
\crefname{conjecture}{Conjecture}{Theorems}
\Crefname{conjecture}{Conjecture}{Theorems}
\crefname{prop}{Proposition}{Propositions}
\Crefname{prop}{Proposition}{Propositions}
\crefname{cor}{Corollary}{Corollaries}
\Crefname{cor}{Corollary}{Corollaries}
\crefname{lem}{Lemma}{Lemmas}
\Crefname{lem}{Lemma}{Lemmas}
\theoremstyle{definition}
\crefname{defn}{definition}{definitions}
\Crefname{defn}{Definition}{Definitions}
\crefname{conj}{Conjecture}{Conjectures}
\Crefname{conj}{Conjecture}{Conjectures}
\crefname{remark}{Remark}{Remarks}
\Crefname{remark}{Remark}{Remarks}
\crefname{rmk}{Remark}{Remarks}
\Crefname{rmk}{Remark}{Remarks}
\crefname{example}{Example}{Examples}
\Crefname{example}{Example}{Examples}
\crefname{align}{}{}
\Crefname{align}{}{}
\crefname{equation}{}{}
\Crefname{equation}{}{}

\def\eqref#1{{\normalfont(\ref{#1})}}

\usepackage{authblk}

\author[1]{Hao Hu\footnote{School of Mathematical and Statistical Sciences, Clemson University, Clemson, USA; Email: \url{hhu2@clemson.edu}; Research supported by the Air Force Office of Scientific Research under award number FA9550-23-1-0508.}}

\begin{document}
\title{The Maximum Singularity Degree for Linear and Semidefinite Programming}
\date{\today}
\maketitle

\begin{abstract}
Facial reduction (FR) is an important tool in linear and semidefinite programming, providing both algorithmic and theoretical insights into these problems. The maximum length of an FR sequence for a convex set is referred to as the maximum singularity degree (MSD). The MSD gives a choice-robust worst-case bound on the number of nontrivial FR steps. It also yields a sufficient rank for a low-rank formulation of the SDP exposing-vector search, while upper bounds on the MSD can serve as proof devices for bounding the singularity degree in structured problem classes. These concrete roles motivate our study of its fundamental properties.

In this work, we show that if an FR sequence has the longest length, then it satisfies a certain minimal property. For linear programming (LP), we prove that every minimal FR sequence forms a basis of a fixed vector space. This yields a direct characterization of the longest FR sequences.  To study the MSD for semidefinite programming (SDP), we provide several useful tools including simplification and upper-bounding techniques. By leveraging these tools and the characterization for LP problems, we prove that finding a longest  FR sequence for SDP problems is NP-hard. This complexity result highlights a striking difference between the shortest and the longest FR sequences for SDP problems.
\end{abstract}

\noindent\textbf{Key Words:}
semidefinite program, linear program, facial reduction, singularity degree, maximum singularity degree, exposing vector
\medskip

\noindent\textbf{Mathematics Subject Classification:} 65K05, 90C27

\newpage

\section{Introduction}
Let $L \cap K$ be a convex set defined as the intersection of an affine subspace $L$ and a closed convex cone $K$. The description of $L \cap K$ is important in efficiently solving optimization problems constrained by $L \cap K$. A poor description of $L \cap K$ can lead to numerical issues, resulting in unreliable outcomes. To address these challenges, Borwein and Wolkowicz introduced a theoretical framework known as the \emph{Facial Reduction Algorithm} (FRA) in their works \cite{borwein1981regularizing,borwein1981facial,borwein1981characterization}. FRA addresses the loss of strict feasibility by reformulating the feasible set over its minimal face. Since its inception, the FRA has significantly improved the computational efficiency and stability of solving numerous \emph{semidefinite programming} (SDP) problems, while also providing deeper insights into the geometric properties of the associated convex sets.

From a computational perspective, the effectiveness of the FRA in numerically solving SDP problems was first demonstrated through its application to the quadratic assignment problem in the study by Zhao et al. \cite{zhao1998semidefinite}. Recent advancements in first-order methods, especially when combined with the FRA, have enhanced the ability to solve SDP relaxations for large-scale problems, as shown in \cite{zheng2017fast,gaar2023strong,oliveira2018admm,hu2020solving,wiegele2022sdp}.

Informally, an FR sequence records the successive steps in which FRA restricts the cone to smaller faces until strict feasibility is restored. From a theoretical perspective, the lengths of FR sequences offer valuable insights into various important questions. Sturm \cite{sturm2000error} introduces a parameter called the \emph{singularity degree} (SD) which is the minimum length among all FR sequences for $L \cap K$. The SD provides an important error bound for $L \cap K$, and has been instrumental in deriving significant theoretical results for different convex cones, as discussed in \cite{drusvyatskiy2017note,pataki2021exponential,lourencco2021amenable,lindstrom2020error,tanigawa2017singularity}. 
The maximum length among all FR sequences for $L\cap K$ is called the \emph{maximum singularity degree} (MSD). The direct literature on the MSD is limited. The strengthened Barvinok--Pataki bound of Im and Wolkowicz \cite{im2021strengthened} is stated in terms of the usual singularity degree; in subsequent work, Im and Wolkowicz \cite[Footnote~4]{im2023revisiting} observe that an MSD-based version can also be obtained. The latter work studies the MSD directly for \emph{linear programming (LP)} and shows that long FR sequences can negatively affect numerical algorithms. Among the references discussed here, it is the only work that studies the MSD directly. The results in \cite{ito2017bound,lourencco2018facial} concern lengths of chains of faces and provide related geometric background, but they do not study the MSD itself. Additionally, we provide a new application of the MSD in Section \ref{sec_app} to further motivate our study.

The SD and MSD also describe two complementary extremes of FRA. If an FR sequence has length $\ell$, then
\[
\sd(L\cap K)\leq \ell\leq \msd(L\cap K).
\]
For SDP, choosing maximum-rank exposing vectors at successive steps is a standard way to attain the lower extreme, but finding such vectors may require additional computational effort. In contrast, any implementation that produces a valid nontrivial exposing vector at each iteration terminates within $\msd(L\cap K)$ steps, regardless of the choices made. Thus, the MSD gives a choice-robust worst-case bound for FRA; this interpretation is structural and does not presume that the MSD is efficiently computable from the problem data.

The MSD can also serve as a proof device without being computed for individual instances. In subsequent work on equality-generated SDP--RLT relaxations of binary programs with $n\geq2$ variables \cite[Section~5.3, Theorem~5.1 and Corollary~5.1]{hu2026rlt}, the worst-case singularity degree is shown to be exactly $\lfloor n/2\rfloor$. The proof uses an inductive argument that restricts FR sequences to faces and bounds the number of nontrivial steps that remain. Since the resulting sequences need not be shortest, the induction requires an upper bound on the MSD of the restricted problems. Thus, the MSD provides the uniform control needed to establish the sharp result.


%
 
The main contribution of this work is to address the problem of constructing the longest FR sequences for linear and semidefinite programming. While the construction of the shortest FR sequences has been thoroughly investigated, extending this analysis to the longest FR sequences provides a natural and important direction for further investigation.

The paper is organized as follows. We provide some preliminary materials in Section \ref{sec_prel}. We also discuss a novel application of the MSD in Section \ref{sec_app}, which serves to further motivate our study. In Section \ref{minFR}, we demonstrate that an FR sequence having the longest length necessarily satisfies a minimal property, and we establish sufficient conditions for minimality in LP and SDP. In Section \ref{sec_msd_lp}, we use the span identity from Corollary~\ref{minseqspan} to prove that every minimal FR sequence for an LP problem forms a basis of a fixed vector space. Thus all minimal FR sequences have the same length, which proves that an FR sequence is longest if and only if it is minimal. In Section \ref{sec_SDP}, we provide counterexamples showing that, for SDP problems, an FR sequence satisfying the minimal property does not necessarily have the longest length. Furthermore, we prove that finding a longest FR sequence for SDP problems is NP-hard, a result that highlights a striking difference between these two parameters, the MSD and the SD.

\textbf{Notation:} Let $\mathbb{R}^n$ denote the $n$-dimensional real space, and let $\mathbb{R}^n_+$ represent the $n$-dimensional real space with nonnegative entries. Consider a finite set $\mathcal{N}$ with $n$ elements. Let $w \in \mathbb{R}^n$. When the entries of $w$ are indexed by the elements of $\mathcal{N}$, we may alternatively write $w \in \mathbb{R}^\mathcal{N}$. For any subset $S \subseteq \mathcal{N}$, let $w(S)$ denote the subvector of $w$ consisting of entries indexed by the elements in $S$. For example, if $\mathcal{N} = \{1, \ldots, n\}$ follows the standard labeling and $n \geq 3$, and if $S = \{1, 2, 3\}$, then $w(S)$ is the vector in $\mathbb{R}^3$ containing the first three entries of $w$. 

Let $F$ and $G$ be subsets of $\Rn$. The \emph{support} of $F$, denoted by $\supp(F)$, is the subset of $\{1,\ldots,n\}$ such that $i \in \supp(F)$ if and only if there exist some $x \in F$ such that the $i$-th entry of $x$ is nonzero. The set difference $F \setminus G$ is defined by $F \setminus G := \{ x \mid x \in F \text{ and } x \notin G\}$.

Let $\mathbb{S}^n$ denote the set of $n \times n$ \emph{symmetric matrices}, and let $\mathbb{S}^n_+$ represent the set of $n \times n$ \emph{positive semidefinite matrices}. Given $X, Y \in \mathbb{S}^n$, the \emph{trace inner product} between $X$ and $Y$ is defined as $\langle X, Y \rangle = \mathrm{tr}(XY)$, where $\tr(\cdot)$ denotes the trace of a matrix.  If the rows and columns of \(X \in \mathbb{S}^n\) are indexed by the elements of a finite set \(\mathcal{N}\), we may alternatively write \(X \in \mathbb{S}^{\mathcal{N}}\). For any subset \(S \subseteq \mathcal{N}\), the principal submatrix of \(X\) corresponding to \(S\) is denoted by \(X(S, S)\).

For any set $L$, the orthogonal complement of $L$ is denoted by $L^{\perp}$. For singleton sets, we simplify the notation by writing  $w^{\perp}$ instead of $\{w\}^{\perp}$.

\section{Preliminaries and Motivation}\label{sec_prel}
\subsection{Facial reduction algorithm (FRA)}
Here and throughout, $\ri(\cdot)$ and $\conv(\cdot)$ denote the relative interior and convex hull, respectively.
Let $K$ be a nonempty closed convex cone in a finite dimensional Euclidean space. The \emph{dual cone} of $K$ is $K^{*} = \{ y \mid \langle y , x \rangle \geq 0, \forall x \in K\}$. We say a nonempty convex subcone $F\subseteq K$ is a \emph{face} of $K$, denoted by $F \unlhd K$, if $x,y \in K$ and $x+y \in F$ imply that $x,y \in F$. Throughout, we consider only nonempty faces. If $F$ is a face of $K$ such that $F$ is nonempty and $F \neq K$, then we say $F$ is a \emph{proper face}. A face $F$ of $K$ is called \emph{exposed} if it is of the form $F = K \cap v^{\perp}$ for some $v \in K^{*}$. The vector $v$ is then called an \emph{exposing vector}. We say $K$ is \emph{exposed} if all of its faces are exposed. 
The conjugate face of $F$ is  $F^{\triangle} := K^{*} \cap F^{\perp}$. For any $x \in \ri F$, we have $F^{\triangle} = K^{*}\cap x^{\perp}$; see, e.g., \cite{tunccel2012strong}. For any $S \subseteq K$, the smallest face of $K$ containing $S$ is denoted by $\face(S,K)$. Note that $\face(S,K) = \face(\conv(S),K)$. In particular, for any $x \in \ri(\conv(S))$, we have $\face(S,K) =\face(x,K)$; see, e.g., Proposition~2.2.5 in \cite{cheung2013preprocessing}.

Let $L$ be an affine subspace such that $L \cap K \neq \emptyset$. We say \emph{Slater's condition} holds for $L \cap K$ if it contains a feasible solution in the relative interior of $K$, i.e., $L \cap \ri K \neq \emptyset$.  The smallest face of $K$ containing $L \cap K$ is called the \emph{minimal face} of $L \cap K$, see  \cite{liu2018exact}. To find the minimal face of $L \cap K$, FRA exploits the following theorem of the alternative
$$L \cap \ri K = \emptyset \quad  \Leftrightarrow \quad  L^{\perp}  \cap (K^{*} \setminus K^{\perp}) \neq \emptyset.$$
We describe FRA applied to $L \cap K$ below, and its proof of convergence can also be found in some recent works, see \cite{liu2018exact,pataki2013strong,waki2013facial}.

\begin{algorithm}[H]
	\caption{Facial Reduction Algorithm (FRA)}\label{alg1}
	\begin{algorithmic}[1]
		\State \textbf{Initialization:} Let $F_0 = K$, $i = 1$.
		\While{there exists $w_i \in L^\perp \cap (F_{i-1}^{*} \setminus F_{i-1}^\perp)$}
		\State Set $F_i \leftarrow F_{i-1} \cap w_i^\perp$.
		\State Set $i \leftarrow i + 1$.
		\EndWhile
	\end{algorithmic}
\end{algorithm}

If $L \cap \ri K = \emptyset$, then FRA applied to $L \cap K$ generates a sequence of faces $(F_{0},\ldots,F_{d})$ satisfying
$$F_{0} \supsetneq F_{1} \supsetneq \cdots \supsetneq F_{d},$$
where  $F_{d}$ is the minimal face of $L \cap K$. It also generates a sequence of exposing vectors $(w_{1},\ldots,w_{d})$ satisfying
$$w_{i} \in  L^\perp \cap (F_{i-1}^{*} \setminus F_{i-1}^\perp) \text{ for } i =1,\ldots,d.$$
We call $(w_{1},\ldots,w_{d})$ an FR sequence for $L \cap K$, corresponding to the sequence of faces $(F_{0},\ldots,F_{d})$. The positive integer $d$ is called the length of the FR sequence. If $L \cap \ri K \neq \emptyset$, then FRA terminates immediately, and we set $d = 0$.

The \emph{singularity degree of $L \cap K$} is defined as the minimum length among all FR sequences for $L \cap K$. For linear and semidefinite programming problems, a shortest FR sequence can be obtained by selecting an exposing vector
\[
w_i\in \ri\!\left(L^\perp\cap F_{i-1}^*\right)\setminus F_{i-1}^\perp
\]
at each FR step. For $K = \Snp$ with $n\geq2$, the singularity degree can be any integer between $0$ and $n-1$. For $K = \Rnp$, the singularity degree is either $0$ or $1$.

The \emph{maximum singularity degree (MSD)} of $L \cap K$, denoted by $\msd(L \cap K)$, is the maximum length among all FR sequences for $L \cap K$. In contrast to the singularity degree, numerous fundamental questions about the MSD remain open, such as how to identify the longest FR sequences for both LP and SDP problems. This paper aims to address these questions.

Throughout, we make the following assumptions to avoid trivial cases in our analysis. Unless otherwise stated or ambiguity arises, the symbols $K$ and $L$ will consistently denote a closed convex cone and an affine subspace, respectively, without repeated definitions.

\begin{assumption}\label{mainass}
Let $L$ be an affine subspace and $K$ a closed convex cone such that:
	\begin{enumerate}
		\item $K$ and $K^{*}$ are both nonempty and exposed.
		\item $L \cap K \neq \emptyset$.
	\end{enumerate}
\end{assumption}

In addition, we rely on the following well-known results, which serve as essential tools for analyzing the relationships between convex sets and their faces.
\begin{lem}\label{conjugate}
	\begin{enumerate}
		\item\label{conjugate1} Let $\emptyset \neq S \subseteq K$ and $F\unlhd K$. Then
		\[
		F=\face(S,K)
		\quad\Longleftrightarrow\quad
		\ri(\conv(S))\subseteq\ri(F).
		\]
		\item Let $F \unlhd K$. Then $F$ is an exposed face if and only if $F^{\triangle\triangle} = F$.
		\item Let $F,G \unlhd K$ be exposed. If $F \subsetneq G$, then $F^{\triangle} \supsetneq G^{\triangle}$.\label{conjugate3}
		\item Let $F,G \unlhd K$. If $F \subsetneq G$, then $\dim F < \dim G$. \label{conjugate4}
	\end{enumerate}
\end{lem}
\begin{proof}
\begin{enumerate}
	\item This is a well-known result; see, e.g., Proposition 2.2.5 in \cite{cheung2013preprocessing}.
	\item This result appears in Proposition 3.1, Item 2, of \cite{tunccel2012strong}.
	\item We have $F^{\triangle} \supseteq G^{\triangle}$ from the definition of the conjugate face. If $F^{\triangle} = G^{\triangle}$, then we have $F = (F^{\triangle})^{\triangle} = (G^{\triangle})^{\triangle} = G$. This contradiction proves the result.
	\item See Corollary 5.5 in \cite{brondsted2012introduction}. \qedhere
\end{enumerate}
\end{proof}

\subsection{FRA for LP}
Let $K = \Rnp$. Define the affine subspace $L:=  \{x \in \Rn \mid	Ax = b \}$ for some $A \in \R^{m \times n}$ and $b \in \R^{m}$. Then $L \cap \Rnp$ is a polyhedron, and we obtain an LP problem.
The entries of $x$ are indexed by $\{1,\ldots,n\}$. Recall that, for any subset $S \subseteq \{1,\ldots,n\}$, $x(S) \in \R^{S}$ denotes the subvector of $x$ indexed by $S$.
A set $F$ is a nonempty face of the nonnegative orthant $\R_{+}^{n}$ if and only if there exists a subset $S \subseteq \{1,\ldots,n\}$ such that 
\begin{equation}\label{rnface}
	F = \left\{ x \in \Rnp \mid x(S) = 0\right\}.
\end{equation}
In the FRA, the set $L^{\perp} \cap F^{*}$ is also a polyhedron. To see this, note that the orthogonal complement of $L$ and the dual cone of the face $F$ are given by 
\begin{equation}\label{lp_dual}
	\begin{array}{rll}
		L^{\perp} &=& \left\{ A^{T}y \in \Rn \mid b^{T}y = 0 \right\},\\
		F^{*} &=& \left\{ w \in \Rn \mid  w(\{1,\ldots,n\}\setminus S) \geq 0 \right\}.
	\end{array}
\end{equation}

\subsection{FRA for SDP} \label{frasdp}
Let $K = \Snp$. Let $A_{1},\ldots,A_{m} \in \Sn$ and $b \in \Rm$ be given. Consider the following SDP problem
\begin{equation}\label{sdp_def}
	L \cap \Snp \text{ where } L := \{ X \in \Sn \mid \langle A_{i}, X \rangle = b_{i} \text{ for } i =1,\ldots,m\}.
\end{equation}
A set $F$ is a nonempty face of $\Snp$ if and only if there exists a linear subspace $\mathcal{V}$ of $\Rn$ such that
\begin{equation}\label{sdpface}
	F=\left\{X \in \Snp  \mid \text{range}(X) \subseteq \mathcal{V}\right\}.
\end{equation}
See, e.g., \cite[Section~2]{pataki2013strong} for this characterization.
Here, $\text{range}(X)$ denotes the range space of $X$. In the FRA, the set $L^{\perp} \cap F^{*}$ also defines an SDP problem. Let $V$ be any matrix such that $\text{range}(V) = \mathcal{V}$. Then
\begin{equation*}
	\begin{array}{rll}
		L^{\perp} &=& \left\{ \sum_{i=1}^{m} A_{i}y_{i} \in \Sn \mid b^{T}y = 0 \right\},\\
		F^{*} &=& \left\{W \in \Sn \mid V^{T}WV  \text{ is positive semidefinite}\right\}.
	\end{array}
\end{equation*}

A special case is when all the matrices in the face $F$ have a block-diagonal structure, possibly after some reordering. Assume the rows and columns of $X$ are indexed by the elements of a finite set $\mathcal{N}$. Let $S \subseteq \mathcal{N}$. The set of positive semidefinite matrices whose rows and columns  corresponding to elements in $\mathcal{N}\setminus S$ are zero is a face, i.e.,
\begin{equation}\label{blk}
F = \left\{X \in \Snp  \mid  X(\mathcal{N}\setminus S,\mathcal{N}\setminus S) = 0\right\}.
\end{equation}
Note that $X(\mathcal{N}\setminus S, \mathcal{N}\setminus S)$ is the principal submatrix of $X$ corresponding to $\mathcal{N}\setminus S$. For example, let \(\mathcal{N} = \{1, \ldots, n\}\) be the standard labeling of the rows and columns from \(1\) to \(n\). Let \(S = \{1, \ldots, r\}\). Then \(F\) can be written as
\[
\begin{array}{rll}
	F &=& \left\{ X \in \Snp \mid X = 
	\begin{bmatrix}
		R & 0 \\
		0 & 0
	\end{bmatrix} \text{ with } R \in \Srp \right\}.
\end{array}
\]
In the above case, the dual cone $F^{*}$ admits a simple characterization: $W \in F^{*}$ if and only if the leading $r \times r$ principal submatrix of $W$ is positive semidefinite.

In our main theorem, we construct an SDP problem with a special structure, ensuring that all faces in any FR sequences exhibit a block-diagonal structure. This allows us to leverage the simple structures of $F$ and $F^{*}$.

\subsection{An application of the MSD}\label{sec_app}
In this section, we provide a new application of the MSD to further motivate our study. A central question in the implementation of FRA is how to find an element in \(L^{\perp} \cap (K^{*}\setminus K^{\perp})\) efficiently. For instance, when \(K = \Snp\), finding an element in \(L^{\perp} \cap \Snp\) becomes an SDP problem, making it a challenging task. This question is very important, and thus, many special FRAs have been developed to achieve efficient implementations in practice; see \cite{permenter2018partial,zhu2019sieve,hu2023affine,friberg2016facial,cheung2013preprocessing}.

Here, we provide a discussion on a special FRA based on the well-known low-rank approach for SDP problems, highlighting how the MSD offers valuable insights into the behavior of this algorithm as an application. The low-rank approach for SDP problems, also known as the Burer-Monteiro SDP method, was introduced in \cite{burer2003nonlinear}, and it has demonstrated significant practical success and has received considerable attention over the past decade.

As mentioned in Section \ref{frasdp}, finding an element in $L^{\perp} \cap (\Snp \backslash \{0\})$ is equivalent to searching for a vector $y = (y_{i})_{i=1}^{m} \in \Rm$ such that
\begin{equation}\label{exp}
	\sum_{i=1}^{m} A_{i}y_{i} \in \Snp, \tr\left(\sum_{i=1}^{m} A_{i}y_{i}\right) = 1 \text{ and } b^{T}y = 0.
\end{equation}
To overcome the difficulties in solving \eqref{exp}, we introduce a new matrix variable $V \in \R^{n \times r}$ for some positive integer $r \leq n$ and consider the following non-linear system in the variable $(y,V) \in \Rm \times \mathbb{R}^{n \times r}$,
\begin{equation}\label{nl}
	\sum_{i=1}^{m} A_{i}y_{i} = VV^{T}, \tr\left(\sum_{i=1}^{m} A_{i}y_{i}\right) = 1 \text{ and } b^{T}y = 0.
\end{equation}
In this factorized formulation, the positive semidefinite constraint in
\eqref{exp} is replaced by smooth equations involving the low-rank factor
$V$. When $r$ is much smaller than $n$, storing and optimizing over $V$
may be computationally attractive compared with introducing an explicit
matrix variable in $\Sn$. This motivates using nonlinear least-squares
methods, such as the Gauss--Newton or Levenberg--Marquardt algorithm, to
search for a solution. Whenever a solution $(y^{*},V^{*})$ is found,
$y^{*}$ is feasible for \eqref{exp} and provides the exposing vector
$\sum_{i=1}^{m}A_i y_i^{*}$ for $L\cap\Snp$.

However, if \(r\) is too small, the non-linear system \eqref{nl} may fail to capture any exposing vectors for \(L \cap \Snp\). Specifically, let \(r^{*}\) denote the smallest rank of exposing vectors from \eqref{exp} for \(L \cap \Snp\). If \(r < r^{*}\), then there is no feasible solution for \eqref{nl}, and consequently, it fails to detect any exposing vectors.  The behavior of  the Burer-Monteiro SDP method has been investigated from many different aspects; see, e.g., \cite{burer2005local,o2022burer,waldspurger2020rank}. The MSD yields the following general upper bound on $r^{*}$, and hence a sufficient rank parameter for \eqref{nl}.

\begin{prop}\label{prop:msd_rank_bound}
Assume that $L\cap\Snp$ does not satisfy Slater's condition. Let
\[
\rho:=\max\{\rank(X)\mid X\in L\cap\Snp\}
\qquad\text{and}\qquad
d:=\msd(L\cap\Snp).
\]
Then
\[
r^{*}\leq n-\rho-d+1.
\]
Consequently, the non-linear system \eqref{nl} has a nonzero solution whenever
$r\geq n-\rho-d+1$. In particular, if $L\cap\Snp$ contains a nonzero matrix,
then $r^{*}\leq n-d$.
\end{prop}

\begin{proof}
Choose a longest FR sequence $(w_1,\ldots,w_d)$ with corresponding
faces $(F_0,\ldots,F_d)$, and let $\mathcal V_i$ be the subspace
associated with $F_i$ in \eqref{sdpface}. Then
$\dim\mathcal V_0=n$ and $\dim\mathcal V_d=\rho$, since the
minimal face contains a feasible matrix in its relative interior.
The first step reduces the subspace dimension by $\rank(w_1)$,
and each subsequent step reduces it by at least one. Thus
\[
n-\rho\geq \rank(w_1)+(d-1).
\]
After positive scaling, $w_1$ satisfies \eqref{exp}, so
$r^*\leq\rank(w_1)\leq n-\rho-d+1$.
The claim for \eqref{nl} follows from PSD factorization.
Finally, a nonzero feasible matrix implies $\rho\geq1$, giving
$r^*\leq n-d$. \qedhere
\end{proof}

Writing $s:=\sd(L\cap\Snp)$, we have $s\leq d$, and hence
$r^*\leq n-\rho-d+1\leq n-\rho-s+1$.
Thus, the MSD-based sufficient rank bound is at least as tight as the
analogous SD-based bound, and is strictly tighter when $d>s$.

The proposition is an existence result: it guarantees that the rank parameter
does not exclude all exposing vectors, but it does not guarantee that a
nonlinear algorithm will recover such a vector reliably or stably.

The following problem is a well-known example for clarifying the numerical issues of ill-conditioned SDP problems, see \cite{tunccel2016polyhedral,sremac2021error,permenter2018partial}. Define the affine subspace
\[
L = \left\{ X \in \Sn \mid X_{11} = 1, \, X_{22} = 0, \, X_{k+1,k+1} = X_{1,k} \text{ for } k = 2, \ldots, n-1 \right\}.
\]
The singularity degree of \(L \cap \Snp\), with \(L\) defined above, is \(n-1\), representing the worst possible scenario. However, the MSD of $L \cap \Snp$ is also \(n-1\). The maximum rank of a feasible matrix is $\rho=1$, so Proposition~\ref{prop:msd_rank_bound} gives $r^{*}\leq1$ and hence $r^{*}=1$. Thus, the non-linear system \eqref{nl} contains a rank-one exposing vector. If a nonzero solution of the corresponding low-rank system can be found at every iteration, the resulting exposing vectors reach the minimal face of \(L \cap \Snp\) in $n-1$ FR steps. This conclusion concerns the existence of rank-one exposing vectors; numerically recovering them, especially for ill-conditioned SDPs, is a separate issue and is not guaranteed by the rank bound.

\section{The minimal FR sequences}\label{minFR}

Let $K$ be a closed convex cone and $L$ an affine subspace.  Let $F \unlhd K$. We say $w$ is \emph{minimal} for $L \cap F$ if $w \in  L^{\perp} \cap (F^{*} \setminus F^{\perp})$ and there does not exist $u \in L^{\perp} \cap (F^{*} \setminus F^{\perp})$ such that
\begin{equation}\label{mindef}
	F \cap w^{\perp} \subsetneq F \cap u^{\perp}.
\end{equation}
Let $f = (w_{1},\ldots,w_{d})$ be an FR sequence for $L \cap K$, and $(F_{0},\ldots,F_{d})$ the corresponding sequence of faces. We call $f$ \emph{minimal} if $w_{i}$ is minimal for $L \cap F_{i-1}$ for all $i=1,\ldots,d$. Our first result is that any longest FR sequence is necessarily minimal. 

\begin{thm}
	\label{lem_split_sdp}
	If $f$ is a longest FR sequence for $L \cap K$, then $f$ is minimal.
\end{thm}
\begin{proof}
	Let $f = (w_{1},\ldots,w_{d})$ be an FR sequence for $L \cap K$, and $(F_{0},\ldots,F_{d})$ the corresponding sequence of faces. Assume that $w_i$ is not minimal for $L \cap F_{i-1}$. By definition, there exists a $u \in L^\perp \cap (F_{i-1}^* \setminus F_{i-1}^\perp)$ such that
\begin{equation}\label{min1}
	F_{i} \subsetneq G \subsetneq F_{i-1} \text{ with } G:= F_{i-1}\cap u^{\perp}.
\end{equation}
	Since taking duals reverses the inclusion order, the above inclusion implies that
	\[
	w_i \in F_{i-1}^* \subseteq G^*.
	\]
In addition, the two inclusions in \eqref{min1} also imply that
	$$G\cap w_{i}^{\perp} = G \cap (F_{i-1} \cap w_i^\perp)  =  G \cap F_{i} = F_{i}.$$
	This means $w_{i} \in L^{\perp} \cap (G^{*} \setminus G^{\perp})$. Thus $(w_1, w_2, \ldots, w_{i-1}, u, w_i, \ldots, w_d) $ is an FR sequence for $L\cap K$ of length $d + 1$. Thus, $f$ is not a longest FR sequence for $L\cap K$. \qedhere
\end{proof}

For LP and SDP problems, we provide two simple sufficient conditions ensuring that an exposing vector $w$ is minimal for $L \cap F$. This result is required in Lemma \ref{allone} for $\Rnp$.
\begin{lem}\label{dim1}
Let $F \unlhd K$ and $w \in L^{\perp} \cap (F^{*}\setminus F^{\perp})$ be an exposing vector for $L \cap F$. 
\begin{enumerate}
	\item If $K = \Rnp$ and $\dim (F) - \dim (F \cap w^{\perp}) = 1$, then $w$ is minimal.
	\item If $K = \Snp$ and $\max\{ \rank(X) \mid X \in F \} - \max\{ \rank(X) \mid X \in F \cap w^{\perp} \} = 1$, then $w$ is minimal.
\end{enumerate}
\end{lem}
\begin{proof}
	The result follows directly from  Item \ref{conjugate4} in Lemma \ref{conjugate}, along with the characterization of the faces of $\Rnp$ and $\Snp$ in \eqref{rnface} and \eqref{sdpface}. \qedhere
\end{proof}

\section{Maximum Singularity Degree for Linear Programming}\label{sec_msd_lp}

In this section, we prove that an FR sequence for a polyhedron is minimal
if and only if it is longest. We first establish a relative support property
of minimal exposing vectors for the nonnegative orthant and then derive a
span identity for minimal FR sequences. Together, these results show that
every minimal FR sequence forms a basis of $L^\perp\cap F_{\min}^\perp$,
where $F_{\min}$ is the minimal face containing the feasible set.
The resulting LP characterization will be used in the proof of
Theorem~\ref{main_hard}.

\begin{lem}\label{sdpminunique}
Let $F\unlhd\Rnp$ be nonempty, set $I:=\supp(F)$, and suppose that
\[
	w\in L^\perp\cap(F^*\setminus F^\perp)
\]
is minimal for $L\cap F$. If $z\in L^\perp\setminus F^\perp$ satisfies
\begin{equation}\label{eq:relative-support}
	\supp(z)\cap I\subseteq\supp(w)\cap I,
\end{equation}
then there exists a nonzero scalar $\alpha\in\R$ such that
\[
	z-\alpha w\in F^\perp.
\]
\end{lem}
\begin{proof}
	Suppose that there is no scalar $\alpha$ such that
	$z-\alpha w\in F^\perp$.
	Define a face $G$ of $F^*$ by
	\[
		G:=\{v\in F^*: \supp(v)\cap I\subseteq\supp(w)\cap I\}.
	\]
	Then $w\in\ri(G)$. Moreover, condition~\eqref{eq:relative-support} implies that
	$z\in\operatorname{span}(G)$. Since $z$ is not proportional to $w$
	modulo $F^\perp$, moving from $w$ in one of the directions $z$ or $-z$
	until reaching the relative boundary of $G$ produces, for some nonzero
		scalar $\beta$, the vector
	\[
		u:=w-\beta z
	\]
		such that $u\in L^\perp\cap(F^*\setminus F^\perp)$ and
	\[
		\supp(u)\cap I\subsetneq\supp(w)\cap I.
	\]
	Consequently,
	\[
		F\cap w^\perp\subsetneq F\cap u^\perp,
	\]
	contradicting the minimality of $w$ for $L\cap F$. Therefore
	$z-\alpha w\in F^\perp$ for some scalar $\alpha$. This scalar is
	nonzero because $z\notin F^\perp$. \qedhere
\end{proof}

\begin{cor}
\label{minseqspan}
Let $f$ be a minimal FR sequence for $L\cap\Rnp$, written as
$f=(w_1,\ldots,w_d)$, and let $(F_0,\ldots,F_d)$ be its corresponding
faces. Then, for every $i\in\{0,\ldots,d\}$,
\begin{equation}\label{eq:minseqspan}
	L^\perp\cap F_i^\perp
	=\operatorname{span}\{w_1,\ldots,w_i\},
\end{equation}
where the span is understood to be $\{0\}$ when $i=0$.
\end{cor}
\begin{proof}
	Set $H_i:=L^\perp\cap F_i^\perp$ for $i=0,\ldots,d$. Since
	$F_0=\Rnp$, we have $H_0=\{0\}$. We claim that
	\[
		H_i=H_{i-1}+\operatorname{span}\{w_i\},
		\qquad i=1,\ldots,d.
	\]
	Indeed, $F_i\subseteq F_{i-1}$ implies
	$H_{i-1}\subseteq H_i$. Moreover,
	$F_i=F_{i-1}\cap w_i^\perp\subseteq w_i^\perp$ and
	$w_i\in L^\perp$, so $w_i\in H_i$. Therefore
	\[
		H_{i-1}+\operatorname{span}\{w_i\}\subseteq H_i.
	\]
	For the reverse inclusion, let
	$z\in H_i\setminus H_{i-1}$ and set $I:=\supp(F_{i-1})$. Because
	$F_i=F_{i-1}\cap w_i^\perp$ and $F_{i-1}$ is a face of the
	nonnegative orthant, the indices removed from the support at the $i$-th step are
	exactly $\supp(w_i)\cap I$. Thus,
	\[
		I\setminus\supp(F_i)=\supp(w_i)\cap I.
	\]
	Moreover, $z\in H_i\subseteq F_i^\perp$, so
	$\supp(z)\cap\supp(F_i)=\emptyset$. Consequently,
	\[
		\supp(z)\cap I
		\subseteq I\setminus\supp(F_i)
		=\supp(w_i)\cap I.
	\]
	Furthermore, $z\in H_i\subseteq L^\perp$, whereas
	$z\notin H_{i-1}=L^\perp\cap F_{i-1}^\perp$. Hence
	$z\notin F_{i-1}^\perp$, and therefore
	$z\in L^\perp\setminus F_{i-1}^\perp$. Since $w_i$ is minimal for
	$L\cap F_{i-1}$, the support containment above and
	Lemma~\ref{sdpminunique} yield a nonzero scalar $\alpha$ such that
	$z-\alpha w_i\in F_{i-1}^\perp$. It also belongs to $L^\perp$, and
	hence to $H_{i-1}$. Thus,
	\[
		z=(z-\alpha w_i)+\alpha w_i
		\in H_{i-1}+\operatorname{span}\{w_i\}.
	\]
	This proves the reverse inclusion and hence the claimed equality. Iterating it from
	$H_0=\{0\}$ gives \eqref{eq:minseqspan} for every
	$i\in\{0,\ldots,d\}$. \qedhere
\end{proof}

\begin{thm}\label{mainthm}
Let $f$ be an FR sequence for the polyhedron $L\cap\Rnp$.
Then $f$ is minimal if and only if it is a longest FR sequence.
\end{thm}

\begin{proof}
Every longest FR sequence is minimal by Theorem~\ref{lem_split_sdp}.
Conversely, let $f=(w_1,\ldots,w_d)$ be a minimal FR sequence, with
corresponding faces $(F_0,\ldots,F_d)$, where $F_0=\Rnp$ and
$F_d=F_{\min}$. Corollary~\ref{minseqspan} gives
\[
L^\perp\cap F_{\min}^\perp
=\operatorname{span}\{w_1,\ldots,w_d\}.
\]
For every $i=1,\ldots,d$ and $j<i$, we have
$F_{i-1}\subseteq F_j\subseteq w_j^\perp$. Hence
\[
\operatorname{span}\{w_1,\ldots,w_{i-1}\}
\subseteq F_{i-1}^\perp.
\]
Since the $i$-th FR step is strict, $w_i\notin F_{i-1}^\perp$.
Thus each $w_i$ lies outside the span of its predecessors, so
$w_1,\ldots,w_d$ are linearly independent. Consequently,
\[
d=\dim\bigl(L^\perp\cap F_{\min}^\perp\bigr).
\]
This also holds for the empty sequence when $d=0$, by Corollary~\ref{minseqspan}.
The right-hand side is independent of the chosen sequence, so every minimal
FR sequence has the same length. Since a longest FR sequence is minimal,
$f$ is longest. \qedhere
\end{proof}

The span identity in Corollary~\ref{minseqspan} need not hold for SDP.
Although the inclusion
\[
\operatorname{span}\{w_1,\ldots,w_i\}\subseteq L^\perp\cap F_i^\perp
\]
holds for every FR sequence over any closed
convex cone, the reverse inclusion can fail even for a minimal sequence. Indeed, let
$K=\mathbb S_+^2$ and $L=\operatorname{span}\{E_{22}\}$, where
$E_{ij}$ denotes the standard matrix unit. The only nonzero matrices
in $L^\perp\cap K$ are positive multiples of $E_{11}$, so
$w_1=E_{11}$ defines a minimal FR sequence ending at
$F_1=\operatorname{cone}\{E_{22}\}$. However,
\[
L^\perp\cap F_1^\perp
=\operatorname{span}\{E_{11},E_{12}+E_{21}\}
\supsetneq\operatorname{span}\{E_{11}\}.
\]
In this example, the sequence is both minimal and longest.

\section{Maximum Singularity Degree for Semidefinite Programming}\label{sec_SDP}
In contrast to the LP characterization in Theorem~\ref{mainthm}, a minimal FR sequence for an SDP problem need not be longest. In fact, we show that even if an exposing vector of minimum rank is chosen at each FR step, this may not yield a longest FR sequence. We prove that the complexity of finding a longest FR sequence is NP-hard. This observation highlights a fundamental difference between the shortest and the longest FR sequences.

\subsection{Examples}
\label{sec_minFR}
Consider the SDP problem \eqref{sdp_def}  defined by the following data matrices,
\begin{equation}\label{notminex}
A_{1} := \begin{bmatrix}
	1 & 0 & 0\\
	0 & 0 & 0\\
	0 & 0 & 0
\end{bmatrix}, A_{2}:=\begin{bmatrix}
	-1 & 1 & 0\\
	1 & 1 & 0\\
	0 & 0 & 0
\end{bmatrix}, A_{3}:=\begin{bmatrix}
	0 & 0 & 0\\
	0 & 1 & 0\\
	0 & 0 & 1
\end{bmatrix} \text{ and } b := \begin{bmatrix}
	0\\
	0\\
	0\\
\end{bmatrix}\in \R^{3}.
\end{equation}
The only feasible solution is zero. The set of exposing vectors $L^{\perp} \cap \mathbb{S}_{+}^{3}$ is any positive semidefinite matrix of the following form
$$ \begin{bmatrix}
y_{1}-y_{2} & y_{2} & 0\\
y_{2} & y_{2} + y_{3} & 0\\
0 & 0 & y_{3}
\end{bmatrix} \succeq 0 \text{ for some } y_{1},y_{2},y_{3} \in \R.$$
It is straightforward to verify that the sequence of exposing vectors $(A_{3},A_{1})$ is a minimal FR sequence, and it has length $2$. However, it is not a longest FR sequence, as the FR sequence $(A_{1},A_{2},A_{3})$ is minimal by Lemma \ref{dim1} and it has length $3$.

In fact, selecting an exposing vector of minimum rank at each FR step does not always yield a longest FR sequence. We illustrate this with an example.  Consider the SDP problem with data matrices,
\begin{equation}\label{SDPex2}
A_{1} := \begin{bmatrix}
1 & 0 & 0 & 0 & 0\\
0 & 1 & 0 & 0 & 0\\
0 & 0 & 1 & 0 & 0 \\
0 & 0 & 0 & 0 & 0 \\
0 & 0 & 0 & 0 & 0 \\
\end{bmatrix}, A_{2}:= \begin{bmatrix}
0 & 0 & 0 & 0 & 0\\
0 & 0 & 0 & 0 & 0\\
0 & 0 & 0 & 1 & 0 \\
0 & 0 & 1 & 1 & 0 \\
0 & 0 & 0 & 0 & 0 \\
\end{bmatrix}, A_{3}:= \begin{bmatrix}                                                                                                                                                                                                                                                                                                                                                                                                                                                
0 & 0 & 0 & 0 & 0\\
0 & 0 & 0 & 0 & 0\\
0 & 0 & 0 & 0 & 0 \\
0 & 0 & 0 & 1 & 0 \\
0 & 0 & 0 & 0 & 1 \\
\end{bmatrix} \text{ and } b := \begin{bmatrix}
0\\
0\\
0\\
\end{bmatrix}\in \R^{3}.
\end{equation}
In the first FR step, all possible exposing vectors up to positive scaling and their ranks can be listed as follows
$$\begin{array}{c|c}
\text{exposing vector} & \text{ rank } \\ \hline
A_{3} & 2\\
A_{1}+ \alpha A_{2} \text{ for } \alpha \in \{0,1\} & 3\\
A_{1}+ \alpha A_{2} \text{ for } \alpha \in (0,1) & 4\\

A_{1}+ \alpha A_{2}+ \beta A_{3} \succeq 0 \text{ for } \alpha \in \{\frac{1-\sqrt{1+4\beta}}{2},\frac{1+\sqrt{1+4\beta}}{2}\}, \beta > 0&  4\\ 
A_{1}+ \alpha A_{2}+ \beta A_{3} \succeq 0 \text{ for } \alpha \in (\frac{1-\sqrt{1+4\beta}}{2},\frac{1+\sqrt{1+4\beta}}{2}) , \beta > 0&  5\\ 
\end{array}$$
Based on the minimum rank rule, $A_{3}$ should be used as the exposing vector in the first FR step. In the second FR step, the leading $3\times3$ principal submatrix of every valid exposing vector is a positive multiple of $I_3$. Thus, all such vectors expose the same face, and we may choose $A_1$, which has minimum rank $3$. This yields an FR sequence $(A_{3},A_{1})$ of length $2$. However, $(A_{1},A_{2},A_{3})$ is a longer FR sequence, and $A_{1}$ does not have the minimum rank.

\subsection{Main ideas and tools}\label{sec_ub}
In this section, we introduce the key tools for proving the main result in Theorem \ref{main_hard}. A central idea in the proof is to establish an upper bound for the MSD of a special SDP problem. To achieve this, we first reduce the given SDP problem to an LP problem with the same MSD. The MSD of this LP problem can then be upper bounded by the MSD of simpler LP problems. By applying Theorem~\ref{mainthm}, we derive a formula for the MSD of these simpler LP problems, leading to the desired upper bound. 

To accomplish this, we present three key tools in this section.
 
\subsubsection{Simplification}\label{sec_sim}
If the data matrices of an SDP problem have special structures, the problem can often be simplified, making its MSD easier to compute. The following results follow directly and are provided without proof.

\begin{lem}\label{simp1}
	Let $L=\{X \in \Sn \mid \langle A_{i}, X \rangle = b_{i} \text{ for } i = 1,\ldots,m \}$ for some $A_{i} \in \Sn$ $(i=1,\ldots,m)$ and $b\in \Rm$. 
	\begin{enumerate}
		\item Assume $F\unlhd \Snp$ has a block-diagonal structure given by
		\begin{equation}\label{simp1eq1}
				F =\left\{X \in \Snp  \mid X = \begin{bmatrix}
				R & 0\\
				0 & 0
			\end{bmatrix} \text{ with } R \in 	\mathbb{S}_{+}^{k+1} \right\}.
		\end{equation}

		 Let $\tilde{A}_{i}$ be the $(k+1)$-th leading principal submatrix of $A_{i}$. Define $$\tilde{L} := \{ \tilde{X} \in \mathbb{S}^{k+1} \mid \langle \tilde{A}_{i}, \tilde{X} \rangle = b_{i} \text{ for } i = 1,\ldots,m \}.$$ Then $\msd(L \cap F) = \msd(\tilde{L} \cap \mathbb{S}_{+}^{k+1})$.\label{simp_i1}
		
		\item {Assume} $\tilde{A}_{1},\ldots,\tilde{A}_{m}$ are diagonal matrices. Let $\tilde{a}_{i} = \text{diag}(\tilde A_{i}) \in \mathbb{R}^{k+1}$. Define $$\tilde{H} := \{ \tilde{x} \in \mathbb{R}^{k+1} \mid  \tilde{a}_{i}^{T}\tilde{x} = b_{i} \text{ for } i = 1,\ldots,m\}.$$ Then 
		$\msd(\tilde{L} \cap \mathbb{S}_{+}^{k+1}) = \msd(\tilde{H} \cap \R_{+}^{k+1}).$\label{simp_i2}
		
		\item Assume that the last entry of $\tilde{a}_{i}$ is zero, i.e., $\tilde{a}_{i} = \begin{bmatrix}
			a_{i}\\
			0
		\end{bmatrix} \in \R^{k+1}$ where $a_{i} \in \R^{k}$, for all $i = 1,\ldots,m$. Define
		$$H := \{ x \in \R^{k} \mid  a_{i}^{T}x = b_{i} \text{ for } i = 1,\ldots,m\}.$$ 
		Then  $\msd(\tilde{H} \cap \R_{+}^{k+1}) = \msd(H \cap \R_{+}^{k}).$\label{simp_i3}
	\end{enumerate}
\end{lem}

\subsubsection{An upper bound}\label{msdub2}
In this section, we establish an upper bound for the MSD. Intuitively, this bound allows us to decompose the problem into smaller, more manageable subproblems, whose MSD can be derived analytically. While the upper bound is presented in a general setting, for the proof of Theorem \ref{main_hard}, we only require its special case when $K$ is the nonnegative orthant.

\begin{lem}\label{key}
	Assume that $K$ is pointed, i.e., $K\cap(-K)=\{0\}$, and that $L \cap K = \{0\}$. Let $G^{1},\ldots,G^{p}$ be faces of $K$ such that \begin{equation}\label{ubass}
		\face(\bigcup_{j=1}^{p}(G^{j}\cap F),F) = F \quad\forall F \unlhd K.
	\end{equation}
	Then $$\msd(L \cap K) \leq \sum_{i=1}^{p} \msd(L \cap G^{i})\label{simp1_i1}.$$
\end{lem}
\begin{proof}
	Let $f=(w_{1},\ldots,w_{d})$ be a longest FR sequence for
	$L\cap K$, of length $d=\msd(L\cap K)$, and let
	$(F_{0},\ldots,F_{d})$ be the corresponding sequence of faces.
	For any nonempty $G\unlhd K$, define a sequence of faces as
	follows. Since $L\cap K=\{0\}$ and $G\subseteq K$, we have
	$L\cap G=\{0\}$. Let $G_{0} = G$. For $i = 1,\ldots,d$, define $G_{i} = G_{i-1} \cap w_{i}^{\perp}$. Since $G_{i} \unlhd F_{i}$, we have $F_{i}^{*} \subseteq G_{i}^{*}$. Thus, $w_{i} \in L^{\perp} \cap (F_{i-1}^{*} \setminus F_{i-1}^{\perp})$ implies that $w_{i} \in L^{\perp} \cap G_{i-1}^{*}$. This implies that $(G_{0},\ldots,G_{d})$ is a sequence of faces such that $G_{i} \unlhd G_{i-1}$ for $i = 1,\ldots,d$. However, it is possible that $w_{i} \in G_{i-1}^{\perp}$ for some $i$, and in this case, we have $G_{i} = G_{i-1}$. Moreover,
	$G_i=G\cap F_i$ for every $i$. Since $K$ is pointed, $\{0\}$ is a face
	of $K$. Hence $L\cap K=\{0\}$ implies that its minimal face is
	$F_d=\{0\}$, and therefore $G_d=\{0\}$. Deleting the non-strict steps yields an FR
	sequence for $L\cap G$. Its length is the number of strict
	containments and is consequently bounded above by $\msd(L\cap G)$, i.e.,
	$$|\{ i \in \{1,\ldots,d\} : G_{i} \subsetneq G_{i-1} \}| \leq \msd(L \cap G).$$
	For each $j \in \{1,\ldots,p\}$, we construct the sequence of faces $(G_{0}^{j},\ldots,G_{d}^{j})$ for $L \cap G^{j}$ as described above. Note that $G_{0}^{j} = G^{j}$. For every $i =1,\ldots,d$, we have $$G_{i}^{j} = G^{j} \cap w_{1}^{\perp} \cap \cdots \cap w_{i}^{\perp} = G^{j} \cap K \cap w_{1}^{\perp} \cap \cdots \cap w_{i}^{\perp} = G^{j} \cap F_{i}.$$
	Thus, by \eqref{ubass}, we have
	\begin{equation}\label{key_eq1}
		\face(\bigcup_{j=1}^{p}G_{i}^{j},F_{i}) = \face(\bigcup_{j=1}^{p}(G^{j} \cap F_{i}),F_{i}) = F_{i}.
	\end{equation}
	Next we show that the above condition implies that, for any $i \in \{1,\ldots,d\}$, there exists at least one index $j \in \{1,\ldots,p\}$ such that $w_{i} \notin (G_{i-1}^{j})^{\perp}$ and thus the containment below is strict
	 $$G_{i}^{j} = G_{i-1}^{j} \cap w_{i}^{\perp} \subsetneq G_{i-1}^{j}.$$ 
	Let $x \in \ri \left(\conv (\bigcup_{j=1}^{p}G_{i-1}^{j})\right) \subseteq \ri (F_{i-1})$, see \eqref{key_eq1} and Item \ref{conjugate1} in Lemma \ref{conjugate}. Then $x = \sum_{j=1}^{p} \lambda_{j}v_{j}$ for some nonzero $\lambda \geq 0$ and $v_{j} \in G_{i-1}^{j}$. As $w_{i} \in F_{i-1}^{*} \setminus F_{i-1}^{\perp}$, we have 
	$$0 < \langle x,w_{i} \rangle = \sum_{j=1}^{p} \lambda_{j} \langle v_{j} , w_{i}\rangle.$$
	As $G_{i-1}^{j} \unlhd F_{i-1}$, we have $F_{i-1}^{*} \subseteq (G_{i-1}^{j})^{*}$ and thus $\langle v_{j} , w_{i}\rangle \geq 0$. This means $\langle v_{j},w_{i} \rangle > 0$ for at least one index $j$. This shows that $w_{i} \notin (G_{i-1}^{j})^{\perp}$ and thus $G_{i}^{j} \subsetneq G_{i-1}^{j}.$ This yields the desired inequality as
	\[
	\msd(L \cap K) \leq \sum_{j=1}^{p}|\{ i \in \{1,\ldots,d\} : G_{i}^{j} \subsetneq G_{i-1}^{j} \}|  \leq \sum_{j=1}^{p} \msd(L \cap G^{j}).
	\qedhere
	\]
\end{proof}

In the proof of Theorem \ref{main_hard}, we need to apply Lemma \ref{key} to a polyhedron. We show that it suffices to find faces $G^{1},\ldots,G^{p}$ of $\Rnp$ such that $\face(\bigcup_{j=1}^{p}G^{j},\Rnp) = \Rnp$. Then the condition \eqref{ubass} holds automatically.

\begin{lem}\label{ub_lp}
	Assume $G^{1},\ldots,G^{p}$ are faces of $\Rnp$ such that 
	\begin{equation}\label{ub_lpass}
	\face(\bigcup_{j=1}^{p}G^{j},\Rnp) = \Rnp.
	\end{equation}
	 Then
	$\face(\bigcup_{j=1}^{p}(G^{j}\cap F),F) = F$ for any $F \unlhd \Rnp$.
\end{lem}
\begin{proof}
	By the facial structure of $\Rnp$ in \eqref{rnface}, there exist subsets $T,S^{1},\ldots,S^{p} \subseteq \{1,\ldots,n\}$ such that
	$$\begin{array}{rll}
		F &=& \{ x \in \Rnp \mid x(T) = 0\},\\
		G^{i} &=& \{ x \in \Rnp \mid x(S^{i}) = 0\} \text{ for } i=1,\ldots,p.\\
	\end{array}$$
	We have $G^{i} \cap F = \{ x \in \Rnp \mid x(S^{i} \cup T) = 0\}.$
	The assumption \eqref{ub_lpass} is equivalent to $\cap_{i=1}^{p} S^{i} = \emptyset$. This implies that $\cap_{i=1}^{p} (S^{i}\cup T) = T$ and thus $\face(\bigcup_{j=1}^{p}(G^{j}\cap F),F) = F$. \qedhere
\end{proof}

We note that Lemma \ref{ub_lp} does not hold for the cone of positive semidefinite matrices. Let $K = \mathbb{S}_{+}^{3}$. Recall that the nonempty faces of $\Snp$ are characterized by the linear subspaces of $\Rn$, see \eqref{sdpface}. Let $G^{1},G^{2}$ be faces of $\mathbb{S}_{+}^{3}$ associated with the following linear subspaces, respectively.
$$ \spa \begin{bmatrix}
	1 & 0\\
	0 & 1\\
	0 & 0\\
\end{bmatrix} \text{ and }\spa \begin{bmatrix}
0 \\
0 \\
1 \\
\end{bmatrix}.$$
Clearly, we have $\face(G^{1}\cup G^{2},\mathbb{S}_{+}^{3}) = \mathbb{S}_{+}^{3}$.
Let $W$ be the all-ones matrix. Then $W \in K^{*} = \mathbb{S}_{+}^{3}$. The exposed face $F = \mathbb{S}_{+}^{3} \cap W^{\perp}$ is associated with the linear subspace
$$\spa \begin{bmatrix}
	1 & 1 \\
	-1 & 0\\
	0 & -1
\end{bmatrix}.$$
However, the faces $G^{1} \cap F$ and $G^{2} \cap F$ are associated with the following linear subspaces, respectively.
$$\spa  \begin{bmatrix}
	1\\
	-1\\
	0
\end{bmatrix} \text{ and } \{0\}.$$
Clearly, the smallest face of $\mathbb{S}_{+}^{3}$ containing both $G^{1} \cap F$ and $G^{2} \cap F$ is just $G^{1} \cap F$, but $G^{1} \cap F \subsetneq F$. This shows  $\face((G^{1}\cap F)\cup(G^{2}\cap F),F) \subsetneq  F$.

\subsubsection{The MSD of a special set}\label{msdspecial}
In this section, we derive the MSD of a special polyhedron.	Let $M = (M_{ij}) \in \{0,1\}^{p \times q}$ be a binary matrix with $p$ rows and $q$ columns. Let $\mathcal{E} = \{ (i,j) \mid M_{ij} = 1\}$ denote the set of indices corresponding to nonzero entries in $M$. We define the polyhedron $H \cap \R_{+}^{\mathcal{E}}$, where the associated affine subspace $H$ is given by
	\begin{equation}\label{exM1}
	H:=	\{ x \in \R^{\mathcal{E}} \mid a_{i}^{T}x = 0 \text{ for } i = 1,\ldots,p+q\}
	\end{equation} for some binary vectors $a_{i} \in \{0,1\}^{\mathcal{E}}$. Note that the entries in the vectors $a_{1},\ldots,a_{p+q}$ and $x$ are indexed by the elements in $\mathcal{E}$. For $i = 1,\ldots,p$, the ones in $a_{i} \in \R^{\mathcal{E}}$ are exactly the entries associated with the ones in the $i$-th row of $M$. For $j = 1,\ldots,q$, the ones in $a_{p+j} \in \R^{\mathcal{E}}$ are exactly the entries associated with the ones in the $j$-th column of $M$. For example, let 
	$$M = \begin{bmatrix}
		1 & 0 & 1\\
		1 & 1 & 0
	\end{bmatrix} \in \R^{2 \times 3}.$$
	Then $p=2$, $q=3$, and
	\[
	\mathcal{E}=\{(1,1),(1,3),(2,1),(2,2)\}.
	\]
	The vectors $a_{1},\ldots,a_{5}$ and $x$ are indexed by $\mathcal E$ and
	can be written as follows:
	\[
	\begin{aligned}
	x&=\begin{bmatrix}x_{(1,1)}\\x_{(1,3)}\\x_{(2,1)}\\x_{(2,2)}\end{bmatrix}
	\in\R^{\mathcal E},
	& a_{1}&=\begin{bmatrix}1\\1\\0\\0\end{bmatrix},
	& a_{2}&=\begin{bmatrix}0\\0\\1\\1\end{bmatrix},\\[0.5em]
	a_{3}&=\begin{bmatrix}1\\0\\1\\0\end{bmatrix},
	& a_{4}&=\begin{bmatrix}0\\0\\0\\1\end{bmatrix},
	& a_{5}&=\begin{bmatrix}0\\1\\0\\0\end{bmatrix}.
	\end{aligned}
	\]
	Note that $a_{1},a_{2}$ correspond to the two rows of $M$, and $a_{3},a_{4},a_{5}$ correspond to the three columns of $M$.

When $M \neq 0$, the only feasible solution to $H \cap \R_{+}^{\mathcal{E}}$ is the zero vector. In general, $a_{i}$ may be a zero vector if it corresponds to a zero row or zero column in $M$; and we retain such vectors only for notational convenience. We show that if $M$ is an all-ones matrix, then the MSD of the corresponding polyhedron $H \cap \R_{+}^{\mathcal{E}}$ can be derived analytically. This result relies on the characterization of the longest FR sequences for polyhedra in Theorem~\ref{mainthm}.

\begin{lem}\label{allone}
	Let $M \in \{0,1\}^{p \times q}$ be an all-ones matrix with $p,q \geq 1$. Let $H \cap \R_{+}^{\mathcal{E}}$ be the polyhedron associated with $M$ as defined in \eqref{exM1}. Then $\msd(H \cap \R_{+}^{\mathcal{E}}) = p+q-1$.
\end{lem}
\begin{proof}
		We construct a minimal FR sequence of length $p + q -1$  for $H \cap \R_{+}^{\mathcal{E}}$. As $H \cap \R_{+}^{\mathcal{E}}$ is a polyhedron, any minimal FR sequence is also a longest FR sequence (see Theorem~\ref{mainthm}). Note that $a_{1},\ldots,a_{p+q}$ are in the dual cone of any face of $\R_{+}^{\mathcal{E}}$ as they are nonnegative. We distinguish two cases depending on the values of $p$ and $q$.	
		\begin{enumerate}
			\item Assume $p = 1$ or $q = 1$. If $q = 1$, then $(a_{1},\ldots,a_{p})$ is a minimal FR sequence of length $p$. Indeed, one can verify that $a_{i} \in {H^{\perp}} \cap (F_{i-1}^{*}\setminus F_{i-1}^{\perp})$ with $F_{0} =  \R_{+}^{\mathcal{E}}$ and $F_{i} = F_{i-1} \cap a_{i}^{\perp}$. By Lemma \ref{dim1}, $a_{i}$ is minimal as $\dim F_{i-1} - \dim F_{i} = 1$.	 This shows that $\msd(H \cap \R_{+}^{\mathcal{E}}) = p$. Similarly, if $p = 1$, then $\msd(H \cap \R_{+}^{\mathcal{E}}) = q$.
		
			\item Assume $p \geq 2$ and $q \geq 2$. We show that $(a_{1},\ldots,a_{p-1},a_{p+1},\ldots,a_{p+q})$ is a minimal FR sequence. It is straightforward to verify that it is indeed an FR sequence. We show that it is minimal. Let $(F_{0},\ldots,F_{p+q-1})$ be the corresponding sequence of faces. Assume, for contradiction, that $a_{1}$ is not minimal for $F_{0}=\R_{+}^{\mathcal{E}}$. By \eqref{lp_dual} and \eqref{mindef}, there exist  $y_{1},\ldots,y_{p+q}\in \R$ such that $u:= \sum_{i=1}^{p+q} a_{i}y_{i} \in \R_{+}^{\mathcal{E}}\setminus \{0\}$ exposes a proper face of $\R_{+}^{\mathcal{E}}$ and 
		\begin{equation}\label{exS}
			\begin{aligned}
			\emptyset\neq S
			&:=\left\{(i,j)\mid \text{the $(i,j)$-th entry of $u$ is nonzero}\right\}\\
			&\subsetneq\left\{(1,j)\mid j=1,\ldots,q\right\}.
			\end{aligned}
		\end{equation}
			Without loss of generality, we assume that $(1,1) \notin S$. Assume $y_{1}  = \lambda$. Then $y_{p+1} = -\lambda$ as $(1,1) \notin S$. From the containment in \eqref{exS}, we conclude that  $y_{2},\ldots,y_{p} = \lambda$, and therefore, $y_{p+2},\ldots,y_{p+q} = -\lambda$. This leads to $u = 0$, contradicting our assumption.
			
			For $i = 2,\ldots,p-1$, we can apply the above argument in the same way to conclude that $a_{i}$ is minimal for $F_{i-1}$. After that the problem is essentially equivalent to the first case. Thus, for $j = 1,\ldots,q$, we have $a_{p+j}$ is minimal for $F_{p+j-2}$. \qedhere
		\end{enumerate}
\end{proof}

Using the formula for $\msd(H \cap \R_{+}^{\mathcal{E}})$ when $H$ is defined by the all-ones matrix in Lemma \ref{allone}, we provide an upper bound for $\msd(H \cap \R_{+}^{\mathcal{E}})$ when $H$ is defined by an arbitrary binary matrix containing many duplicated columns. We state this upper bound in Lemma \ref{MgeneralUB} in a form that is convenient for the proof of Theorem \ref{main_hard}.

\begin{lem}\label{MgeneralUB}
	Let $p$ and $\tilde q$ be positive integers. For
	$j=1,\ldots,\tilde q$, let $v_j\in\{0,1\}^{p}$ and define
	\[
	M_j:=\begin{bmatrix}v_j&\cdots&v_j\end{bmatrix}
	\in\{0,1\}^{p\times 2\tilde q},
	\]
	where $M_j$ consists of $2\tilde q$ copies of $v_j$. Define
	\[
	M:=\begin{bmatrix}M_1&\cdots&M_{\tilde q}\end{bmatrix}
	\in\{0,1\}^{p\times(2\tilde q^2)}.
	\]
	Let $H\cap\R_{+}^{\mathcal E}$ be the polyhedron associated with $M$ as
	described in \eqref{exM1}.
	\begin{enumerate}
		\item {If $v_{1},\ldots,v_{\tilde q}$ are all nonzero, then
		\[
		\msd(H\cap\R_{+}^{\mathcal E})
		\leq\sum_{j=1}^{\tilde q}(\textbf{1}^{T}v_j+2\tilde q-1).
		\]}
		\item\label{MgeneralUB2} If every $v_j$ contains at most three
		nonzero entries and at least one $v_j$ is zero, then
		\[
		\msd(H\cap\R_{+}^{\mathcal E})
		\leq(2\tilde q+2)(\tilde q-1).
		\]
	\end{enumerate}

\end{lem}
\begin{proof}
	\begin{enumerate}
		\item Recall that $\mathcal{E}$ is the index set of nonzero entries in $M$. Let $\mathcal{E}_{j}$ be the index set of nonzero entries in $M_{j}$, the $j$-th submatrix of $M$. Define the face $G^{j}$ of $\R_{+}^{\mathcal{E}}$  as
		$$G^{j} := \{  x \in \R_{+}^{\mathcal{E}} \mid  x(\mathcal{E}\setminus \mathcal{E}_{j}) = 0 \}.$$
		As $\cap_{j=1}^{\tilde{q}}(\mathcal{E}\setminus \mathcal{E}_{j}) = \emptyset$, it is clear that $\face(\bigcup_{j=1}^{\tilde{q}}G^{j},\R_{+}^{\mathcal{E}}) = \R_{+}^{\mathcal{E}}$. Applying Lemma \ref{key} and Lemma \ref{ub_lp} yields
		$\msd(H \cap \R_{+}^{\mathcal{E}}) \leq \sum_{j=1}^{\tilde{q}} \msd(H \cap G^{j}) .$
		As every column $v_{j}$ is nonzero, $M_{j}$ is an all-ones matrix with possibly some additional rows of zeros. Since the additional rows of zeros do not correspond to any variables, we can apply Lemma \ref{allone} to obtain $\msd(H \cap G^{j})  = \textbf{1}^{T}v_{j} +2\tilde{q} - 1$. This yields the first inequality.
		
		\item If $v_{j}$ is nonzero and contains at most $3$ nonzeros, then we have $\msd(H \cap G^{j}) \leq 2\tilde{q} + 2$. If $v_{j}$ is the zero vector, then $\msd(H \cap G^{j}) = 0$. Since there are at most $\tilde{q}-1$ nonzero columns $v_{j}$, the second inequality follows. \qedhere
	\end{enumerate}
\end{proof}

\subsection{NP-hardness}\label{sec_nphard}
In this section, we show that the decision version of finding a  longest FR sequence for SDP problems is NP-hard.  We will construct a polynomial-time transformation from the well-known NP-complete problem 3SAT. These problems are formally defined below. 

\begin{flushleft}
	\textbf{3SATISFIABILITY (3SAT)}\\
		\textbf{INSTANCE:} A set $U = \{u_{1},\ldots,u_{p}\}$ of Boolean variables, and $C = \{c_{1},\ldots,c_{q}\}$ a collection of clauses on $U$ where each clause $c_{i}$ contains exactly three distinct literals. \\ 
	\textbf{QUESTION:} Is there a truth assignment for $U$ such that all the clauses in $C$ are satisfied?
\end{flushleft}

\begin{flushleft}
	\textbf{MAXIMUM SINGULARITY DEGREE for SDP (MSD-SDP)}\\
	\textbf{INSTANCE:} Given symmetric matrices $A_{1},\ldots,A_{m}\in\mathbb{Q}^{n\times n}$, a vector $b\in\mathbb{Q}^{m}$, and a positive integer $d$, all encoded in binary.\\
	\textbf{QUESTION:} Does $L \cap \Snp$, as defined in \eqref{sdp_def}, admit an FR sequence of length at least $d$?
\end{flushleft}

For our analysis later, we preprocess any given 3SAT instance so that it satisfies some additional assumptions. If a clause	$c_{k} = (u_{i},\bar{u}_{i},u_{j})$ contains both the positive literal $u_{i}$ and the negative literal $\bar{u}_{i}$, then $c_{k}$ is trivially satisfied. Thus we can remove this clause $c_{k}$ from the problem. If the positive literal $u_{i}$ never appears in any clause, then we can assume $u_{i}$ is assigned false, allowing us to remove $u_{i}$ and all clauses containing $\bar{u}_{i}$.  Similarly, if the negative literal $\bar{u}_{i}$ never appears, we can assume $u_{i}$ is true, removing all clauses containing $u_{i}$. We repeat these pure-literal simplifications until neither applies. If all clauses are removed, the original instance is satisfiable, and we output a fixed yes-instance of MSD-SDP. We henceforth consider the case in which at least one clause remains. Since we can implement this preprocessing in polynomial time, we can assume without loss of generality that the given 3SAT instance satisfies the following properties.
\begin{assumption}\label{3satass}
	The 3SAT instance satisfies
	\begin{itemize}
		\item For each variable $u_{i}$, we have $u_{i} \in c_{j}$ and $\bar{u}_{i} \in c_{k}$ for some $j$ and $k$.
		\item Each clause $c_{j}$ contains at most one of $u_{i}$ or $\bar{u}_{i}$, but not both.
	\end{itemize}
\end{assumption}

For any 3SAT instance with $p$ variables and $q$ clauses, we construct an MSD-SDP instance as follows. The order of the matrix variable is $n = 6q + 1$ and the number of constraints is $m=2p+q$.  The rows and columns of the matrix variable $X$ and the data matrices $A_{1},\ldots,A_{m}$ are indexed by elements of the set	
$$\mathcal{N}:=  \mathcal{N}_{1} \cup \mathcal{N}_{2}\cup \{(0,0,0)\},$$
where
$$\mathcal{N}_{k}:= \{(i,j,k) \mid u_{i} \in c_{j} \text{ or } \bar{u}_{i} \in c_{j}, \; i =1,\ldots,p, \;j = 1,\ldots,q \} \text{ for } k = 1,2.$$
Note that $(i,j,k) \in \mathcal{N}_{k}$ for some $k$ if and only if the clause $c_{j}$ contains either the literal $u_{i}$ or the literal $\bar{u}_{i}$. By Assumption \ref{3satass}, it cannot contain both. The third index $k$ is simply used for making two identical copies $\mathcal{N}_{1}$ and $\mathcal{N}_{2}$. The last element $\{(0,0,0)\}$ serves an auxiliary role in the analysis.

Next, we specify the data matrices $A_{1},\ldots,A_{m}$ which can be classified into two different components based on their roles.
\begin{enumerate}
	\item For each \( i \in \{1, \ldots, p\} \), we define two diagonal matrices \( A_{i} \) and \( A_{p+i} \in \Sn \) with entries in $\{0,1\}$. They serve as truth-setting components to enforce a choice between assigning the variable \( u_{i}\in U \) to true or false.
	 Define the subsets
	 \begin{equation}\label{TFV}
\begin{array}{rll}
	\mathcal{T}_{i} &:=& \left\{ (i,j,1) \in \mathcal{N}_{1} \mid  u_{i} \in c_{j} \text{ for some } j \right\},\\
	\mathcal{F}_{i} &:=& \left\{ (i,j,1) \in \mathcal{N}_{1} \mid  \bar{u}_{i} \in c_{j} \text{ for some } j \right\},\\
	\mathcal{V}_{i} &:=& \left\{ (i,j,2)  \in \mathcal{N}_{2} \mid j = 1,\ldots,q\right\}.\\
\end{array}
	 \end{equation}
In all three sets above, the variable $j$ serves as the running index, iterating over clause indices. Conceptually, these sets represent the following: the set $\mathcal{T}_{i}$ can be viewed as the set of clauses containing the positive literal $u_{i}$, while the set $\mathcal{F}_{i}$ can be viewed as the set of clauses containing the negative literal $\bar{u}_{i}$.  The set $\mathcal{V}_{i}$ can be viewed as all clauses containing either the positive literal $u_{i}$ or the negative literal $\bar{u}_{i}$. Finally, the third index in each tuple distinguishes elements in $\mathcal{T}_{i}$ and $\mathcal{F}_{i}$ (which use index $1$) from those in $\mathcal{V}_{i} $ (which use index $2$).


The diagonal entries of $A_{i}$ corresponding to $\mathcal{T}_{i} \cup \mathcal{V}_{i}$, and those of $A_{p+i}$ corresponding to $\mathcal{F}_{i} \cup \mathcal{V}_{i}$, are set to $1$. All other entries are $0$. Note that $A_{i}$ and $A_{p+i}$ are positive semidefinite.


\item For each $j \in \{1,\ldots,q\}$, we define a symmetric matrix $A_{2p+j} \in \Sn$ with binary entries, which serves as a satisfaction testing component for the clause $c_{j}$. Define the subsets
$$\begin{array}{rll}
\mathcal{C}_{j} &:=& \left\{(i,j,1) \in \mathcal{N}_{1} \mid i=1,\ldots,p \right\},\\
\mathcal{D}_{j} &:=& \left\{(i,j,2) \in \mathcal{N}_{2} \mid i=1,\ldots,p \right\}.
\end{array}$$
The sets \(\mathcal{C}_{j}\) and \(\mathcal{D}_{j}\) represent all elements of \(\mathcal{N}_{1}\) and \(\mathcal{N}_{2}\), respectively, that correspond to the clause \(c_{j}\). Note that $|\mathcal{C}_{j}|=|\mathcal{D}_{j}| = 3$.  The nonzero entries of \( A_{2p+j} \) are specified as follows: the diagonal entries corresponding to \(\mathcal{C}_j\) are ones. For the off-diagonal entries, the \((0,0,0)\)-th row and the columns corresponding to \(\mathcal{D}_j\) are set to $1$, with symmetry ensuring the corresponding transpose entries are also $1$. For example, the principal submatrix of $A_{2p+j}$ associated with $\mathcal{D}_{j} \cup \{(0,0,0)\}$ is the matrix
\begin{equation}\label{clause_com}
	\begin{bmatrix}
		0 & 0 & 0 & 1 \\
		0 & 0 & 0 & 1 \\
		0 & 0 & 0 & 1 \\
		1 & 1 & 1 & 0
	\end{bmatrix} \in \mathbb{S}^{4}.
\end{equation}
Note that the principal submatrix \eqref{clause_com} is indefinite and it has rank $2$. The two nonzero eigenvalues are $\pm \sqrt{3}$ corresponding to the eigenvectors $\begin{bmatrix}
	1 & 1 & 1 & \sqrt{3}
\end{bmatrix}^{T}$ and  $\begin{bmatrix}
	 1 & 1 & 1 & -\sqrt{3} 
\end{bmatrix}^{T}$, respectively. 
\end{enumerate}
Let $b = 0 \in \Rm$ be the all-zeros vector of length $m$, and $d := p + q$. This defines an MSD-SDP instance, which can be constructed in polynomial time. 

As usual, the feasible region of the MSD-SDP instance is denoted by $L \cap \Snp$, where $L = \{ X \in \Sn \mid \langle A_{i}, X \rangle = 0 \text{ for } i =1,\ldots,m\}$ is the affine subspace determined by the constructed data matrices. It is not difficult to see that if $X \in L \cap \Snp$, then all rows and columns of $X$ corresponding to $\mathcal{N}_{1}\cup \mathcal{N}_{2}$ are zero, and the $(0,0,0)$-th diagonal entry can be any nonnegative number. Thus, the smallest face of $\Snp$ containing $L \cap \Snp$ has a block-diagonal structure given by
\begin{equation}\label{sf}
	\{ X \in \Snp \mid X(\mathcal{N}_{1}\cup \mathcal{N}_{2},\mathcal{N}_{1}\cup \mathcal{N}_{2}) = 0 \}.
\end{equation}

The sparsity pattern of the matrices in $L^{\perp}$ is important in the subsequent analysis. Let \(W = A_{1}y_{1} + \cdots + A_{m}y_{m} \in L^{\perp}\), where \(y_{1}, \ldots, y_{m} \in \mathbb{R}\). For any \((i,j,1) \in \mathcal{N}_{1}\), the \((i,j,1)\)-th diagonal entry of \(W\) is given by \(y_i + y_{2p+j}\) if the positive literal \(u_i\) appears in clause \(c_j\), and \(y_{p+i} + y_{2p+j}\) if the negative literal \(\bar{u}_i\) appears in clause \(c_j\). Formally,
\begin{equation}\label{Wform1}
	W(\{(i,j,1)\},\{(i,j,1)\}) =
	\begin{cases} 
		y_i + y_{2p+j}, & \text{if } (i,j,1) \in \mathcal{T}_i, \\
		y_{p+i} + y_{2p+j}, & \text{if } (i,j,1) \in \mathcal{F}_i.
	\end{cases}
\end{equation}
For any \((i,j,2) \in \mathcal{N}_{2}\), the \(2\times 2\) principal submatrix of \(W\) corresponding to \(\{(i,j,2), (0,0,0)\}\) is given by
\begin{equation}\label{Wform2}
	\begin{blockarray}{ccc}
		(i,j,2) & (0,0,0) & \\
		\begin{block}{[cc]c}
			y_i + y_{p+i} & y_{2p+j} & (i,j,2) \\
			y_{2p+j} & 0 & (0,0,0) \\
		\end{block}
	\end{blockarray}
\end{equation}

All other entries of \(W\) are zero, except those explicitly defined in \eqref{Wform1} and \eqref{Wform2}.

We clarify the construction and the sparsity pattern using the following concrete example.
\begin{ex}\label{MSDSDPex}
Assume $U=\{u_1,u_2,u_3\}$ and $C=\{c_1\}$, where
$c_1=u_1\vee u_2\vee\bar u_3$. Then
\[
\begin{aligned}
\mathcal N_1&=\{(1,1,1),(2,1,1),(3,1,1)\},\\
\mathcal N_2&=\{(1,1,2),(2,1,2),(3,1,2)\}.
\end{aligned}
\]
For some $y_1,\ldots,y_7\in\R$, the matrix $W\in L^\perp$ has the following
form.
\begin{center}
\resizebox{\textwidth}{!}{$
W = \begin{blockarray}{cccccccc}
	(1,1,1) &  (2,1,1) & (3,1,1) & (1,1,2) & (2,1,2) & (3,1,2) & (0,0,0) \\
	\begin{block}{[ccccccc]c}
		y_{1}+y_{7} & 0 & 0 & 0 & 0 & 0 & 0 & (1,1,1)\\
		0 & y_{2} + y_{7} & 0 & 0 & 0 & 0 &0 & (2,1,1) \\
		0 & 0 & y_{6} + y_{7} & 0 & 0 & 0 & 0 & (3,1,1) \\
		0 & 0 & 0 & y_{1}+y_{4} & 0 & 0 & y_{7} &  (1,1,2) \\
		0 & 0 & 0 & 0 & y_{2}+y_{5} & 0 & y_{7} &  (2,1,2) \\
		0 & 0 & 0 & 0  & 0 & y_{3}+y_{6} & y_{7} &  (3,1,2) \\
		0 & 0 & 0 & y_{7}  & y_{7} & y_{7} & 0 &  (0,0,0) \\
	\end{block}
\end{blockarray}
$}
\end{center}
Note that the 3SAT instance in this example does not satisfy Assumption \ref{3satass}, as it was chosen for illustrative purposes.
\end{ex}

The key idea behind the proof of Theorem \ref{main_hard} is to establish a correspondence between the truth assignments of a given 3SAT instance and the FR sequences for \(L \cap \Snp\) in the constructed MSD-SDP problem. For the forward direction, we provide an informal discussion with a concrete example to clarify the idea. For the backward direction, we show some auxiliary results about the structure of the FR sequences. Finally, we provide the formal proof of both directions in Theorem \ref{main_hard}.

In the forward direction, given a truth assignment that satisfies all clauses in the 3SAT instance, we can construct an FR sequence of length \(p+q\) as follows. First, for the truth-setting components, if \(u_i\) is assigned true, we include \(A_{p+i}\) in the FR sequence; if \(u_i\) is assigned false, we include \(A_{i}\) in the FR sequence. Then, for the satisfaction-testing components, we include \(A_{2p+1}, \ldots, A_{2p+q}\) in the FR sequence. We formally prove the forward direction in Theorem \ref{main_hard}, referencing \eqref{fp2} and \eqref{fp3}. To make it easier for the reader to verify the formal proof, we illustrate the construction using a concrete example here.

To illustrate this construction, consider Example \ref{MSDSDPex}, where \(p = 3\) and \(q = 1\). Suppose a truth assignment is given where \(u_1\) is set to true and \(u_2, u_3\) are set to false. In this case, the clause \(c_{1}\) is satisfied. This truth assignment induces an FR sequence \((A_{4}, A_{2}, A_{3}, A_{7})\) of length \(4\), following the construction process described above. Let \((F_{0}, \ldots, F_{4})\) be the corresponding sequence of faces. These faces exhibit a block-diagonal structure, see \eqref{blk},
$$\begin{array}{rrll}
F_{i} = \{ X \in \Snp \mid X(\mathcal{N} \setminus S_{i}, \mathcal{N}\setminus S_{i}) = 0\},
\end{array}$$
where $S_{0},\ldots,S_{4}$ are subsets of $\mathcal{N}$ given by
\[
\begin{aligned}
	S_{0} &= \{ (0,0,0), (1,1,1), (3,1,1), (3,1,2), (2,1,1), (2,1,2), (1,1,2) \}, \\
	S_{1} &= \{ (0,0,0), (1,1,1), (3,1,1), (3,1,2), (2,1,1), (2,1,2) \}, \\
	S_{2} &= \{ (0,0,0), (1,1,1), (3,1,1), (3,1,2) \}, \\
	S_{3} &= \{ (0,0,0), (1,1,1), (3,1,1) \}, \\
	S_{4} &= \{ (0,0,0) \}.
\end{aligned}
\]

The more challenging direction is to demonstrate that if a 3SAT instance is not satisfiable, then no FR sequence of length \(p+q\) or greater exists. To this end, we first establish that there exists a longest  FR sequence for \(L \cap \Snp\) in the constructed MSD-SDP instance such that it corresponds to a truth assignment of the variables of the 3SAT instance.

We show that the faces in any FR sequence for the constructed MSD-SDP have a block-diagonal structure (see \eqref{blk}). Since the rows and columns of the matrix variable $X$ are indexed by the elements in $\mathcal{N}$, each face in an FR sequence corresponds to a specific subset of $\mathcal{N}$.

\begin{lem}\label{main_lem0}
Let \(f = (W_{1}, \ldots, W_{d})\) be an FR sequence for \(L \cap \Snp\), and let \((F_{0}, \ldots, F_{d})\) denote the corresponding sequence of faces. Define \(S_{0} = \mathcal{N}\). For \(r = 1, \ldots, d\), define subsets as follows:
\begin{equation}\label{Sform}
	S_{r} = \{ (i,j,k) \in S_{r-1} \mid \text{the $(i,j,k)$-th diagonal entry of } W_{r} \text{ is zero}\}.
\end{equation}

For any \(r \in \{0, \ldots, d\}\), the face \(F_{r}\) has a block-diagonal structure given by:
\begin{equation}\label{form1}
	F_{r} = \{ X \in \Snp \mid X(\mathcal{N} \setminus S_{r}, \mathcal{N} \setminus S_{r}) = 0 \}.
\end{equation}
For every \(r\in\{1,\ldots,d\}\), the exposing vector \(W_{r} = A_{1}y_{1} + \cdots + A_{m}y_{m} \in L^{\perp} \cap {\bigl(F_{r-1}^{*}\setminus F_{r-1}^{\perp}\bigr)}\), with \(y_{1}, \ldots, y_{m} \in \mathbb{R}\), satisfies the following property: \footnote{The coefficients \(y_{1}, \ldots, y_{m}\) depend on \(r\), but for readability, we omit explicit dependence on  \(r\) in the notation.}
\begin{equation}\label{yzero}
	(i,j,2) \in S_{r-1} \implies y_{2p+j} = 0.
\end{equation}
In addition, $S_{d} = \{(0,0,0)\}$.
\end{lem}
\begin{proof}
We prove the statement by induction. As \(W_1 \in L^{\perp} \cap \Snp\), it follows from \eqref{Wform2} that \(y_{2p+j} = 0\) for all \(j \in \{1, \ldots, q\}\). Thus, \(W_1\) is a diagonal matrix. Consequently, \(S_1\) and \(F_1\) must be in the forms given by \eqref{Sform} and \eqref{form1}, respectively.

Assume that \(F_{r-1}\) has a block-diagonal structure in \eqref{form1} for some subset \(S_{r-1} \subseteq \mathcal{N}\) satisfying \eqref{Sform}. Then \(W_{r} \in L^{\perp} \cap F_{r-1}^{*}\) implies that the principal submatrix \(W_{r}(S_{r-1}, S_{r-1})\) is positive semidefinite. Since the smallest face of \(\Snp\) containing \(L \cap \Snp\) is given by \eqref{sf}, we have \((0,0,0) \in S_{r-1}\) for all $r$. Thus, for any \((i,j,2) \in S_{r-1}\), the \(2 \times 2\) principal submatrix corresponding to \(\{(i,j,2), (0,0,0)\}\) is positive semidefinite. By \eqref{Wform2}, this implies \(y_{2p+j} = 0\) for any \((i,j,2) \in S_{r-1}\). Hence, \(W_{r}(S_{r-1}, S_{r-1})\) is a diagonal matrix. This ensures that \(F_{r}\) also has the form in \eqref{form1}, with the subset \(S_{r}\) given by $
S_{r} = \{ (i,j,k) \in S_{r-1} \mid \text{the $(i,j,k)$-th diagonal entry of } W_{r} \text{ is zero}\}.$

Since the smallest face of \(\Snp\) containing \(L \cap \Snp\) is given by \eqref{sf}, we conclude that the final subset in the sequence satisfies $S_{d} = \{(0,0,0)\}$. \qedhere
\end{proof}

Using the relation between the subsets of $\mathcal{N}$ and the faces in FR sequences, we establish key properties of FR sequences in the following result. The key observation is that a longest FR sequence can be chosen so that its first $p$ exposing vectors select one of $A_i$ and $A_{p+i}$ for each $i=1,\ldots,p$.

\begin{lem}
\label{main_lem1}
Let $f=(W_1,\ldots,W_d)$ be an FR sequence for $L\cap\Snp$, with
corresponding faces $(F_0,\ldots,F_d)$. Fix $\alpha\in\{1,\ldots,p\}$.
Let $\mathcal T_\alpha$, $\mathcal F_\alpha$, and $\mathcal V_\alpha$ be as
defined in \eqref{TFV}, and let $S_0,\ldots,S_d$ be as in \eqref{Sform}.
Define
\[
l:=\min\left\{r\in\{1,\ldots,d\}:\mathcal{T}_{\alpha}\cup\mathcal{F}_{\alpha}\cup\mathcal{V}_{\alpha}\not\subseteq S_r\right\}.
\]
Then $l$ is well defined and unique, and the following three conditions hold simultaneously:
\begin{enumerate}
	\item $\mathcal{T}_{\alpha} \cup \mathcal{F}_{\alpha} \cup \mathcal{V}_{\alpha} \subseteq S_{l-1}$. \label{main_lem1_i1}
	\item $\mathcal{V}_{\alpha} \cap S_{l} = \emptyset$. \label{main_lem1_i2}
	\item $\mathcal{T}_{\alpha} \cup \mathcal{F}_{\alpha} \not\subseteq S_{l}$. \label{main_lem1_i3}
\end{enumerate}
\end{lem}
\begin{proof}
	By Lemma \ref{main_lem0}, we have $S_0=\mathcal{N}$ and $S_d=\{(0,0,0)\}$. Hence, $\mathcal{T}_{\alpha}\cup\mathcal{F}_{\alpha}\cup\mathcal{V}_{\alpha}$ is contained in $S_0$ but not in $S_d$, so $l$ is well defined. Its definition and the inclusions $S_0\supsetneq\cdots\supsetneq S_d$ imply
	\[
	\mathcal{T}_{\alpha}\cup\mathcal{F}_{\alpha}\cup\mathcal{V}_{\alpha}\subseteq S_{l-1}
	\quad\text{and}\quad
	\mathcal{T}_{\alpha}\cup\mathcal{F}_{\alpha}\cup\mathcal{V}_{\alpha}\not\subseteq S_l.
	\]
	At the $l$-th FR step, write the exposing vector as
	\[
	W_l=\sum_{i=1}^{m}A_i y_i,
	\qquad y_1,\ldots,y_m\in\mathbb R,
	\]
	where $W_l\in L^\perp\cap(F_{l-1}^*\setminus F_{l-1}^\perp)$.
	Since $W_l(S_{l-1},S_{l-1})$ is positive semidefinite and
	$\mathcal V_\alpha\subseteq S_{l-1}$, equation~\eqref{yzero} gives
	\begin{equation}\label{y2pj0}
		\begin{gathered}
		y_{2p+j}=0\quad\text{for every }j\in\{1,\ldots,q\}\\
		\text{such that }(\alpha,j,2)\in\mathcal V_\alpha.
		\end{gathered}
	\end{equation}
	
	By Assumption \ref{3satass}, there exists a clause $c_{\beta}$ containing the positive literal $u_{\alpha}$, i.e., $(\alpha,\beta,1) \in \mathcal{T}_{\alpha}$. Since $W_{l}(S_{l-1}, S_{l-1})$ is positive semidefinite and $\mathcal{T}_{\alpha} \subseteq S_{l-1}$, it follows from \eqref{Wform1} that $y_{\alpha} + y_{2p+\beta} \geq 0$. Substituting $y_{2p+\beta} = 0$ from \eqref{y2pj0}, this implies $y_{\alpha} \geq 0$. Similarly, we can argue that there exists a clause $c_{\gamma}$ containing the negative literal $\bar{u}_{\alpha}$, yielding $y_{p+\alpha} \geq 0$. If $y_{\alpha} = y_{p+\alpha} = 0$, the principal submatrix of $W_{l}$ corresponding to $\mathcal{T}_{\alpha} \cup \mathcal{F}_{\alpha} \cup \mathcal{V}_{\alpha}$ would be zero. Consequently, $\mathcal{T}_{\alpha} \cup \mathcal{F}_{\alpha} \cup \mathcal{V}_{\alpha} \subseteq S_{l}$, contradicting the choice of $l$. Hence, $y_{\alpha} + y_{p+\alpha} > 0$. Consequently, the principal submatrix \(W_{l}(\mathcal{V}_{\alpha}, \mathcal{V}_{\alpha})\) is diagonal with positive diagonal entries $y_{\alpha} + y_{p+\alpha} > 0$. This implies $\mathcal{V}_{\alpha} \cap S_{l} = \emptyset$ by \eqref{Sform}. 	Additionally, if $y_{\alpha} > 0$, then $\mathcal{T}_{\alpha} \cap S_{l} = \emptyset$. Similarly, if $y_{p+\alpha} > 0$, then $\mathcal{F}_{\alpha} \cap S_{l} = \emptyset$. This establishes the three claims. \qedhere
\end{proof}

\begin{lem}\label{main_lem2}
Under the assumptions and notation of Lemma~\ref{main_lem1}, assume further
that $f$ is a longest FR sequence, and let $l$ be the unique index defined in
that lemma. Then:
	\begin{enumerate}
		\item The $l$-th exposing vector can be chosen from
		$\{A_{\alpha},A_{p+\alpha}\}$ up to positive scaling.
		\item With this choice, the reordered sequence
		\[
		(W_{l},W_{1},\ldots,W_{l-1},W_{l+1},\ldots,W_{d})
		\]
		is a longest FR sequence.
	\end{enumerate}
\end{lem}

\begin{proof}
	
	To prove the first claim of this lemma, write $W_l=\sum_{i=1}^{m}A_i y_i$, as in the proof of Lemma~\ref{main_lem1}, and define $\tilde{W}:=A_\alpha y_\alpha+A_{p+\alpha}y_{p+\alpha}$. As established there, $y_\alpha,y_{p+\alpha}\geq0$ and $y_\alpha+y_{p+\alpha}>0$. Since these inequalities hold and $\mathcal{T}_{\alpha} \cup \mathcal{F}_{\alpha} \cup \mathcal{V}_{\alpha} \subseteq S_{l-1}$, it follows that $\tilde{W} \in L^{\perp} \cap (F_{l-1}^{*} \setminus F_{l-1}^{\perp})$. By \eqref{y2pj0} and the sparsity structure of elements in $L^{\perp}$ described in \eqref{Wform1} and \eqref{Wform2}, the principal submatrices of $\tilde{W}$ and $W_{l} - \tilde{W}$ corresponding to $S_{l-1}$ have disjoint nonzero entries. This means $W_{l} -\tilde{W} \in L^{\perp} \cap F_{l-1}^{*}$. Thus, we must have $W_{l} -\tilde{W} \in F_{l-1}^{\perp}$.  If not, the aforementioned nonzero pattern of $\tilde{W}$ and $W_{l} - \tilde{W}$ implies that $W_{l} -\tilde{W} \notin (F_{l-1} \cap \tilde{W}^{\perp})^{\perp}$. Then $(W_{1},\ldots,W_{l-1},\tilde{W},W_{l}-\tilde{W},W_{l+1},\ldots,W_{d})$ is also an FR sequence whose length is $d+1$, which is a contradiction. Consequently,
	\[
	\begin{aligned}
	F_l
	&=F_{l-1}\cap W_l^\perp\\
	&=F_{l-1}\cap(W_l-\tilde W+\tilde W)^\perp\\
	&=F_{l-1}\cap(W_l-\tilde W)^\perp\cap\tilde W^\perp\\
	&=F_{l-1}\cap\tilde W^\perp.
	\end{aligned}
	\]
	To justify the third equality, let $X\in F_{l-1}$. Since
	$\tilde{W},W_l-\tilde{W}\in F_{l-1}^*$, we have
	\[
	\langle X,\tilde{W}\rangle\geq 0
	\quad\text{and}\quad
	\langle X,W_l-\tilde{W}\rangle\geq 0.
	\]
	Consequently,
	\[
	\langle X,W_l-\tilde{W}+\tilde{W}\rangle=0
	\]
		if and only if both inner products above are zero. This proves the third
		equality. The fourth equality follows from $W_{l} - \tilde{W} \in F_{l-1}^{\perp}$. Thus, $W_{l}$ and $\tilde{W}$ expose the same face of $F_{l-1}$. We may therefore replace $W_l$ by $\tilde{W}=A_{\alpha}y_{\alpha}+A_{p+\alpha}y_{p+\alpha}$ without changing $F_l$ or the length of the FR sequence. If $y_{\alpha}$ and $y_{p+\alpha}$ were both positive, then $\mathcal{T}_{\alpha}\cup\mathcal{F}_{\alpha}\subseteq S_{l-1}$ and Assumption~\ref{3satass} would give the FR sequence
		\[
		(W_1,\ldots,W_{l-1},A_\alpha,A_{p+\alpha},W_{l+1},\ldots,W_d)
		\]
		of length $d+1$, a contradiction. Thus exactly one of $y_\alpha$ and
		$y_{p+\alpha}$ is positive. Hence, $\tilde{W}$ is a positive multiple of
		exactly one of $A_{\alpha}$ and $A_{p+\alpha}$. Since positive scaling does
		not change the exposed face, we may normalize this replacement and take
		$W_l\in\{A_{\alpha},A_{p+\alpha}\}$.
	
		For the second claim, let
		\[
		R:=\supp(\operatorname{diag}W_l)=S_{l-1}\setminus S_l.
		\]
		Since $W_l$ is a nonzero diagonal positive semidefinite matrix in
		$L^\perp$, it can be applied first, removing the coordinates in $R$.
		For each $r<l$, we have $R\subseteq S_r$, so the diagonal matrix
		$W_r(S_{r-1},S_{r-1})$ is zero on $R$. Deleting these coordinates
		therefore preserves its positive semidefiniteness and its nonzero
		diagonal entries. Thus $W_r$ remains a valid nontrivial step, leaving
		the index set $S_r\setminus R$. After $W_1,\ldots,W_{l-1}$, this set is
		$S_{l-1}\setminus R=S_l$, so all subsequent steps are unchanged.
		The reordered sequence is therefore valid and has the same length as
		$f$, hence is longest. \qedhere
\end{proof}

We conclude that there exists a longest FR sequence that corresponds to a truth assignment.

\begin{cor}\label{main_lem3}
	There exists at least one longest  FR sequence \((W_{1}, \ldots, W_{d})\) for \(L \cap \Snp\) such that \(W_{i} \in \{A_{i}, A_{p+i}\}\) for  \(i = 1,\ldots,p\).
\end{cor}

\begin{proof}
	Start with any longest FR sequence and apply Lemmas~\ref{main_lem1} and~\ref{main_lem2} successively for $\alpha=p,p-1,\ldots,1$. At each application, the selected exposing vector from $\{A_\alpha,A_{p+\alpha}\}$ is moved to the front, while the vectors selected in the preceding applications shift one position to the right. Each application preserves validity and length. The resulting longest FR sequence therefore satisfies $W_i\in\{A_i,A_{p+i}\}$ for $i=1,\ldots,p$. \qedhere
\end{proof}

We are now ready to prove the main result.

\begin{thm}\label{main_hard}
The MSD-SDP problem is NP-hard.
\end{thm}
\begin{proof}
	Let \(U = \{u_{1}, \ldots, u_{p}\}\) and \(\tilde{C} = \{\tilde{c}_{1}, \ldots, \tilde{c}_{\tilde{q}}\}\) be a 3SAT instance satisfying Assumption \ref{3satass}. Note that \(p \geq 3\). To achieve the desired outcome, we introduce redundancy as follows. Each clause \(\tilde{c}_{i} \in \tilde{C}\) is duplicated to create \(2\tilde{q}\) copies. This results in a new set of clauses given by:
	\begin{equation}\label{main_hard1}
		\{ \underbrace{\tilde{c}_{1}, \ldots, \tilde{c}_{1}}_{2\tilde{q}}, \ldots, \underbrace{\tilde{c}_{\tilde{q}}, \ldots, \tilde{c}_{\tilde{q}}}_{2\tilde{q}} \}.
	\end{equation}
	Define \(q := 2\tilde{q}^{2}\). We consider the new 3SAT instance with the same variable set \(U\) and an expanded clause set \(C = \{c_{1}, \ldots, c_{q}\}\) as defined in \eqref{main_hard1}. For \(j = 1, \ldots, \tilde{q}\), it follows that \(c_{i} = \tilde{c}_{j}\) for any \(i \in \{1 + 2\tilde{q}(j-1), \ldots, 2\tilde{q}j\}\). 

	We then transform this new 3SAT instance into an instance of MSD-SDP following the transformation process outlined earlier in this section. We show that the new 3SAT instance is satisfiable if and only if the constructed MSD-SDP instance has an FR sequence of length at least \(d = p+q\).

Now suppose that the new 3SAT instance has a satisfying assignment.
We construct an FR sequence of length $d$,
\[
f=(W_1,\ldots,W_d),
\]
with corresponding faces $(F_0,\ldots,F_d)$. Let $(S_0,\ldots,S_d)$ be the
associated subsets of $\mathcal N$ defined in \eqref{Sform}.
	
	\begin{itemize}
		\item {The first $p$ steps encode the truth assignment. For
		$i=1,\ldots,p$, choose the exposing vector at step $i$ as follows:}
		\begin{equation}\label{fp2}
			W_{i} = \begin{cases}
				A_{i} & \text{ if $u_{i}$ is false} \\
				A_{p+i} & \text{ if $u_{i}$ is true}. \\
			\end{cases}
		\end{equation}
		By construction, $W_{i} \in L^{\perp} \cap (F_{i-1}^{*}\setminus F_{i-1}^{\perp})$ ensuring that these are valid FR steps. 
		\item The final $q$ FR steps are defined as follows:
				\begin{equation}\label{fp3}
		W_{p+j} = A_{2p+j} \text{ for } j = 1,\ldots,q.
				\end{equation}
		Since the first $p$ FR steps in \eqref{fp2} ensure that $S_{p} \cap \mathcal{N}_{2} = \emptyset$, it follows that the principal submatrix $W_{p+j}(S_{p+j-1},S_{p+j-1})$ is positive semidefinite and $W_{p+j} \in L^{\perp} \cap F_{p+j-1}^{*}$ is an exposing vector for $L \cap F_{p+j-1}$ for $j = 1,\ldots,q$. Thus, it remains to show that $W_{p+j} \notin F_{p+j-1}^{\perp}$ so that it exposes a proper face of $F_{p+j-1}$. By Lemma \ref{main_lem0}, this is equivalent to showing that the principal submatrix of $W_{p+j}$ associated with $S_{p+j-1}$ is nonzero.
	
	Let \(j = 1\). Recall that \(A_{2p+j}\) contains exactly three positive diagonal entries corresponding to the literals in the clause \(c_{j}\). Given a satisfying assignment for the new 3SAT instance, at least one literal in $c_{j}$ is true. Assume the clause $c_{j}$ contains a positive literal \(u_{i}\), and \(u_{i}\) is assigned true. Then \(W_{i} = A_{p+i}\), based on the choice in \eqref{fp2}. Since the \((i,j,1)\)-th diagonal entry is zero in \(A_{p+i}\), it follows that \((i,j,1) \in S_{p+j-1}\). Moreover, the \((i,j,1)\)-th diagonal entry of \(A_{2p+j}\) is one, which implies that $A_{2p+j} \notin F_{p+j-1}^{\perp}$. A similar argument holds if the clause $c_{j}$ contains the negative literal \(\bar{u}_{i}\), and the variable $u_{i}$ is assigned false. This shows that \(W_{p+j}\) yields a valid FR step.
	
	Furthermore, note that \(A_{2p+\beta}\) and \(A_{2p+\gamma}\) do not have any common nonzero entries for distinct \(\beta, \gamma \in \{1, \ldots, q\}\). By repeating this reasoning for each clause, we conclude that \(W_{p+2}, \ldots, W_{p+q}\) are also valid FR steps.
	\end{itemize}
This yields an FR sequence of length $p+q$.

Conversely, assume that the new 3SAT instance is not satisfiable. Let \(f = (W_{1}, \ldots, W_{r})\) be a longest FR sequence for the constructed SDP problem, and let \((F_{0}, \ldots, F_{r})\) denote the corresponding sequence of faces. Let $(S_{0}, \ldots, S_{r})$ be the corresponding sequence of subsets of \(\mathcal{N}\) as defined in  \eqref{Sform}. We prove that the length of $f$ is strictly smaller than \(p + q\). 

By Corollary \ref{main_lem3}, we can assume that
\begin{equation}\label{fp}
	W_{i} \in \{A_{i}, A_{p+i}\} \quad \text{for } i = 1, \ldots, p.
\end{equation}
This induces a truth assignment unambiguously via the relation
\begin{equation}\label{truthass}
	u_{i} = 
	\begin{cases}
		\text{false} & \text{if } W_{i} = A_{i}, \\
		\text{true} & \text{if } W_{i} = A_{p+i}.
	\end{cases}
\end{equation}
The face \(F_{p}\) is given by
\[
F_{p} = \left\{ X \in \Snp \mid {X(\mathcal{N} \setminus S_{p},\mathcal{N} \setminus S_{p}) = 0} \right\},
\]
where
\begin{equation}\label{Edef}
	S_{p} = \left\{ (i,j,1) \in \mathcal{T}_{i} \mid u_{i} \text{ is true} \right\} 
	\cup \left\{ (i,j,1) \in \mathcal{F}_{i} \mid u_{i} \text{ is false} \right\} 
	\cup \{(0,0,0)\}.
\end{equation}
Since \(f\) is a longest FR sequence, its tail
$(W_{p+1},\ldots,W_r)$ must be a longest FR sequence for $L\cap F_p$;
otherwise, replacing it with a longer tail would contradict the maximality of
$f$. Therefore,
\[
\msd(L \cap \Snp) = p + \msd(L \cap F_{p}).
\]

We apply the simplification in Lemma \ref{simp1} to show that \(L \cap F_{p}\) is equivalent to an LP problem, and they have the same maximum singularity degree. Since the face \(F_{p}\) has a block-diagonal structure, we can apply Item \ref{simp_i1} in Lemma \ref{simp1} to simplify the problem by removing the rows and columns corresponding to \(\mathcal{N} \setminus S_{p}\).  Formally, define the data matrices
\[
\tilde{A}_{i} := 
\begin{cases}
	A_{i}(S_{p}, S_{p}) & \text{if } W_{i} = A_{p+i}, \\
	A_{p+i}(S_{p}, S_{p}) & \text{if } W_{i} = A_{i},
\end{cases}
\]
for \(i = 1, \ldots, p\), and
\[
\tilde{A}_{p + j} = A_{2p+j}(S_{p}, S_{p}),
\]
for \(j = 1, \ldots, q\). Let \(\tilde{n} = |S_{p}|\). Define the affine subspace
\[
\tilde{L} := \left\{ \tilde{X} \in \mathbb{S}^{\tilde{n}} \mid \langle \tilde{A}_{i}, \tilde{X} \rangle = 0 \text{ for } i = 1, \ldots, p+q \right\}.
\]
By Item \ref{simp_i1} in Lemma \ref{simp1}, we have $\msd(L \cap F_{p}) = \msd(\tilde{L} \cap \mathbb{S}_{+}^{\tilde{n}})$. Note that we discarded exactly half of the data matrices from \(A_{1}, \ldots, A_{2p}\) in this process. This does not cause any issues. Specifically, if \(W_{i} = A_{i}\), then \(A_{i} \in F_{p}^{\perp}\); if \(W_{i} = A_{p+i}\), then \(A_{p+i} \in F_{p}^{\perp}\). Thus, these matrices can be removed freely.
		
Since \(\tilde{A}_{1}, \ldots, \tilde{A}_{p+q}\) are diagonal matrices, we can further simplify $\tilde{L} \cap \mathbb{S}_{+}^{\tilde{n}}$ by applying Item \ref{simp_i2} in Lemma \ref{simp1}. This yields an equivalent LP problem. Additionally, note that the diagonal entry of \(\tilde{A}_{i}\) associated with \((0,0,0)\) is zero for all $i$. Therefore, Item \ref{simp_i3} in Lemma \ref{simp1} allows us to freely remove the corresponding entry from the LP problem. Define \(\mathcal{E} := S_{p} \setminus \{(0,0,0)\}\), and let $a_{i}$ be the vector of diagonal entries of  $\tilde{A}_{i}(\mathcal{E}, \mathcal{E})$, i.e., \(a_{i} = \text{diag}(\tilde{A}_{i}(\mathcal{E}, \mathcal{E})) \in \mathbb{R}^{\mathcal{E}}\). The final LP problem is \(H \cap \mathbb{R}_{+}^{\mathcal{E}}\), where
\[
H := \left\{ x \in \mathbb{R}^{\mathcal{E}} \mid a_{i}^{T}x = 0 \text{ for } i = 1, \ldots, p+q \right\}.
\]
By Item \ref{simp_i2} and \ref{simp_i3} in Lemma \ref{simp1}, we have $\msd(L \cap F_{p}) = \msd(H \cap \R_{+}^{\mathcal{E}})$. Note that the elements in \(\mathcal{E}\) are uniquely defined by their first two indices $i$ and $j$, as \((i,j,k) \in \mathcal{E}\) implies \(k = 1\). By examining the definition of \(S_{p}\) in \eqref{Edef}, it follows directly that \(H \cap \mathbb{R}_{+}^{\mathcal{E}}\) corresponds to the polyhedron constructed in \eqref{exM1}, which is associated with the following binary matrix \(M = (M_{ij})\in \{0,1\}^{p \times q}\):
\[
M_{i,j} =  
\begin{cases}
	1 & \text{if } u_{i} \text{ is true and } u_{i} \in c_{j}, \\
	1 & \text{if } u_{i} \text{ is false and } \bar{u}_{i} \in c_{j}, \\
	0 & \text{otherwise}.
\end{cases}
\]

Recall that the set of clauses \(C\) contains \(2\tilde{q}\) duplications of each clause in \(\tilde{C}\), as shown in \eqref{main_hard1}. The matrix \(M\) can be written as:
\[
M = \begin{bmatrix}
	M_{1} & \cdots & M_{\tilde{q}}
\end{bmatrix},
\]
where each submatrix \(M_{j} \in \{0,1\}^{p \times 2\tilde{q}}\) corresponds to the \(2\tilde{q}\) duplicates associated with the same clause \(\tilde{c}_{j}\) for \(j = 1, \ldots, \tilde{q}\).

Consider the truth assignment defined in \eqref{truthass}:
\begin{enumerate}
	\item If \(\tilde c_{j}\) is satisfied, then \(M_{j} \neq 0\) and it has at least one row of ones. Additionally, \(M_{j}\) contains at most three rows of ones because each clause contains three literals.
	\item Since the given 3SAT instance is not satisfiable, there exists at least one unsatisfied clause. If \(\tilde c_{j}\) is unsatisfied, then \(M_{j}\) is an all-zeros matrix. Thus, there is at least one submatrix \(M_{j}\) that is equal to zero.
\end{enumerate}
Therefore, the matrix \(M\) satisfies the assumptions in Item \ref{MgeneralUB2} of Lemma \ref{MgeneralUB}. This yields
\begin{equation}\label{npmainineq}
	\msd(H \cap \mathbb{R}_{+}^{\mathcal{E}}) \leq (2\tilde{q} + 2)(\tilde{q} - 1).
\end{equation}
We now obtain an upper bound for the length of \(f\),
\[
\msd(L \cap \Snp) = p + \msd(L \cap F_{p}) \leq p + (2\tilde{q} + 2)(\tilde{q} - 1)  = p + q - 2 < d.
\]
Thus, the constructed MSD-SDP instance does not admit any FR sequences of length \(d\) or more. \qedhere

\end{proof}


%
\section*{Conflict of interest}
The author has no relevant financial or non-financial interests to disclose.




\end{document}